\documentclass[a4paper, 10pt, ]{article}

\usepackage{amssymb,amsmath,enumerate,graphicx,hyperref,color}
 
\usepackage{t1enc}
 
\usepackage[T1]{fontenc}

\usepackage{amsfonts}

\usepackage{color}
\usepackage{graphicx}
\usepackage{rotating}
\usepackage{caption}
\usepackage{fancyhdr}
 
\usepackage{hyperref}
 
\newtheorem{lemma}{Lemma}[section]
\newtheorem{proposition}[lemma]{Proposition}
\newtheorem{theorem}[lemma]{Theorem}

\newtheorem{remark}[lemma]{Remark}
\newtheorem{definition}[lemma]{Definition}
\newtheorem{example}[lemma]{Example}

\newtheorem{assumption}[lemma]{Assumption}

\newcommand{\beq}{\begin{eqnarray}}
\newcommand{\enq}{\end{eqnarray}}
\newcommand{\be}{\begin{eqnarray*}}
\newcommand{\en}{\end{eqnarray*}}
\newcommand{\ben}{\begin{eqnarray*}}
\newcommand{\enn}{\end{eqnarray*}}

\newcommand{\R}{\mathbb R}

\newcommand{\E}{\mathbb E}

\newcommand{\1}{{\mathbf I}}

\newcommand{\F}{\mathbb F}

\let\cal=\mathcal

\newcommand{\Ac}{\mathcal{A}} 
\newcommand{\Fc}{\mathcal{F}}  
\newcommand{\Dc}{\mathcal{D}} 
  
\newcommand{\vr}{{\rm v}}

\newcommand{\eps}{\varepsilon}  
  
\newcommand{\x}{\times}

\def \vp{\varphi}
\def\vs#1{\vspace{#1mm}}

\def\Q{{\mathbb Q}}
\def\vr{{\rm v}}
\def\Ec{{\cal E}}

 \def\vrm{{\rm v}}
 \def\xr{{\rm x}}

  \def\vr{\vrm}
  
  \def\E{{\mathbb E}}
    \def\P{{\mathbb P}}

	\def\Gk{\gamma}
	\def\Bk{B}
	
	\def\Om{\Omega}
	\def\om{\omega}

	\def\Pt{\widetilde \P}
	\def\Xb{\overline X}
	\def\Pcb{\overline\Pc}
	\def\Pct{\widetilde\Pc}
	\def\alphab{\overline \alpha}
	
	\def\Jb{\overline J}
	\def\Ph{\widehat \P}
	\def\Eh{\widehat \E}
	\def\C{\mathbb{C}}
	\def\0{\mathbf{0}}
	\def\Rb{\overline \R}
	\def\gammah{\hat \gamma}
	\def\Gammab{\overline \Gamma}
	\def\Tc{\cal T}
	\def\Pb{\overline \P}

  \def\Pc{{\cal P}}

  \def\R{{\mathbb R}}
    
\def \proof{{\noindent \bf Proof. }}
\def \endproof{\hbox{ }\hfill$\Box$}
\def\dr{{\mathfrak d}}
\def\dr{{d}}
\def\CT{C([0,T])}

\def\DT{D([0,T])}     
\def\TD{[0,T]\x\DT}

 \def\Ninfty#1{\|#1\|}    
 \def\Cb{{\mathbb C}_{\rm r}}      
 \def\Href#1{{\rm (H\ref{#1})}}

\title{Understanding the dual formulation for the hedging of path-dependent options with price impact\thanks{ The authors are grateful to Pierre Cardaliaguet for helpful discussion and suggestions. This work has benefited from the financial support of the Initiative de Recherche ``M\'ethodes non-lin\'eaires pour la gestion des risques financiers'' sponsored by AXA Research Fund.
}
}
\author{
Bruno Bouchard
\footnote{Universit\'e Paris Dauphine - PSL, CNRS, CEREMADE, bouchard@ceremade.dauphine.fr. } 
\and 
Xiaolu Tan
\footnote{Department of Mathematics, The Chinese University of Hong Kong. xiaolu.tan@cuhk.edu.hk}} 

\date{\today}

\begin{document}

\maketitle 

	\begin{abstract}  
		
	\end{abstract}
 
\section{Introduction and notations}
 
	Following the steps of  \cite{AbergelLoeper,LoMi1}, the paper \cite{BoLoZo1} constructed a market model  pertaining for permanent price impact. It is based on a simple linear price impact rule around the origin and a passage to the limit from the situation in which the composition of the portfolio is modified in discrete time. 
	It takes the form 
	\begin{align}
		Y&=Y_{0}+ \int_{0}^{\cdot} a_{t} dW_{t}+\int_{0}^{\cdot} b_{t}dt, \;\;\;\;
		V=V_{0 }+ \int_{0}^{\cdot} Y_{t}dX_{t}+\frac12 \int_{0}^{\cdot} a^{2}_{t} f(X_{t}) dt,
		\label{eq: intro V Y}\\
		X&=X_{0 }+\int_{0}^{\cdot} \mu(X_{t}) dt + \int_{0}^{\cdot} \sigma(X_{t}) dW_{t} + \int_{0}^{\cdot} f(X_{t}) dY_{t}+\int_{0}^{\cdot} a_{t}(\sigma f')(X_{t}) dt\nonumber\\
		&=X_{0 }+\int_{0}^{\cdot} \sigma_{X}^{a_{t}}(X_{t}) dW_{t} + \int_{0}^{\cdot} \mu_{X}^{a_{t},b_{t}}(X_{t})dt,\label{eq: intro X}
	\end{align}
	where 
	$\sigma_{X}^{a_{t}}:=(\sigma+a_{t}f)$ and $\mu_{X}^{a_{t},{b_{t}}}:=( {\mu+}b_{t}f+a_{t}\sigma f' )$, for some functions: $(\mu, \sigma, f): \R \to \R^3$.
	In the above, $Y$ is the number of stocks in the hedging portfolio, controlled by the predictable processes $a$ and $b$, $V$ is the book value of the portfolio (amount of cash plus position in stocks evaluated at their market price) and $X$ is the impacted stock price, given the impact function $f>0$.   
 	In contrast to the un-covered case considered in \cite{BoLoZo1}, the paper \cite{BoLoZo2} studied the problem of hedging a {\sl covered} option of payoff $\Phi(X_{T})$ within this model.
	The terminology {\sl covered} means that the trader asks at time $0$ an amount of cash and the number of stocks he needs to initiate his trading strategy. The stocks are evaluated at their market value to determine the premium. At maturity, the  trader delivers the stocks he has (again evaluated at their market value)  plus an amount in cash to  match the payoff $\Phi(X_{T})$. This avoids having an important impact on the price process when the trader needs to jump to an initial delta at time $0$, and liquidate his current delta at maturity time $T$. Hence, the hedging problem boils down to  finding some constants $(V_{0},Y_{0})$ and processes $(a,b)$ in \eqref{eq: intro V Y}-\eqref{eq: intro X} such that $V_{T}=\Phi(X_{T})$.
	
	\vspace{0.5em}
	
	This can be viewed as a (highly non classical) {\sl second order} coupled forward backward stochastic differential equation (2FBSDE).  The term {\sl second order} comes from the fact that $Y$ itself admits an It\^{o} decomposition as in \cite{cheridito2007second,soner2009dynamic}, with the difference that, in our setting, there is a  coupling    through the quadratic variation of $Y$, which appears both  in the quadratic variation of the forward process $X$  and in the drift of the backward process $V$. 

	\vspace{0.5em}

	The super-hedging problem has been studied in \cite{BoLoZo2} from the stochastic target point of view. 
	The authors provide a viscosity solution characterization  of the super-hedging price function $(t,x)\mapsto v(t,x)$. 
	To understand their result, let us first rewrite the  dynamic of the number of stocks in the form 
$$
Y=Y_{0}+ \int_{0}^{\cdot} \gamma_{t} dX_{t}+\int_{0}^{\cdot} b^{'}_{t}dt
$$ 
in which $\gamma:=a/(\sigma(X)+f(X)a)$ is the {\sl gamma} of the portfolio, and $b^{'}$ is a predictable process.    Assuming  that the hedging strategy consists in tracking the super-hedging price, as in classical complete market models,
then we expect that $V=v(\cdot,X)$. By applying twice It\^{o}'s lemma, this implies $\frac12 a^{2}f(X)=\partial_{t}v(\cdot,X)+\frac12 \sigma_{X}^{a}(X)^{2}\partial^{2}_{xx}v(\cdot,X)$, $Y=\partial_{x}v(\cdot,X)$ and $\gamma=\partial^{2}_{xx}v(\cdot,X)$. Combining the above, we obtain the partial differential equation 
$$
-\partial_{t}v-\frac12 \frac{\sigma^{2}}{1/f-\partial^{2}_{xx}v} \partial^{2}_{xx}v=0
$$
on $[0,T)\x \R$, with terminal condition $v(T,\cdot)=\Phi$. The precise formulation in \cite{BoLoZo2} involves an additional constraint $\partial^{2}_{xx}v\le 1/f$ which ensures that the above PDE remains parabolic, together with the corresponding face-lift of the terminal  condition. It was later observed in \cite{BLSZ18} that, in the case where the   constraint $\partial^{2}_{xx}v< 1/f$ is automatically enforced by the terminal condition $\Phi$ (see Remark \ref{rem: eps conca}  below), then the above PDE is an Hamilton-Jacobi-Bellman equation:
$$
	\inf_{{\rm a}\in \R}  \Big(-\partial_{t}v-\frac12 {{\rm a}^{2}}\partial^{2}_{xx}v+\frac1{2f}({\rm a}-\sigma)^{2} \Big)=0.
$$   
 This in turn leads to the dual formulation 
 	\begin{equation}\label{eq: dual intro}
		v(0,X_{0})
		=
		\sup_{\alpha} \E \Big[\Phi(X^{\alpha}_{T})-\int_{0}^{T}\frac1{2f(X^{\alpha}_{t})} \big(\alpha_{t}-\sigma(X^{\alpha}_{t}) \big)^{2}dt \Big],
		~\mbox{with}~
		X^{\alpha}=X_{0}+\int_{0}^{\cdot} \alpha_{t}dW_{t}.
 	\end{equation}

	The very aim of this paper is to better understand this dual formulation, and to extend it to a more general setting. 
	First, one can observe that it does {not} follow from a super-hedging inequality, as it is usually the case in the literature on option pricing, but rather from a hedging equality. The dual problem is not a natural lower bound, but instead a natural upper bound.  Namely, let us assume that $(V_{0},Y_{0},\hat a,\hat b)$ solves \eqref{eq: intro V Y}-\eqref{eq: intro X} and $V_{T}=\Phi(X_{T})$, and let $\Q^{\hat a,\hat b}$ be a martingale measure for $X$ (assuming that it exists and that $\int_{0}^{\cdot} Y_{s}dX_{s}$ is a martingale). 
	Then,
	\begin{equation*}
		V_{0} 
		~=~ 
		\E^{\Q^{\hat a,\hat b}} \Big[\Phi(X_{T})-\int_{0}^{T} \frac{1}{2}\hat a^{2}_{t}f(X_{t})dt \Big] 
		~\le~ 
		\sup_{a,b}\E^{\Q^{ a, b}} \Big[\Phi(X^{a,b}_{T})-\int_{0}^{T} \frac{1}{2}a^{2}_{t}f(X^{a,b}_{t})dt \Big],
	\end{equation*}
	in which $X^{a,b}=X_{0}+\int_{0}^{\cdot} \sigma_{X}^{a_{t}}(X^{a,b}_{t})dW^{a,b}_{t}$ where  $W^{a,b}$ is a Brownian motion under $\Q^{a,b}$. 
	This is formally equivalent to  
	\begin{equation*}
		V_{0}
		~\le~
		\sup_{a}\E \Big[\Phi(\tilde X^{a}_{T})-\int_{0}^{T} \frac{1}{2}a^{2}_{t}f(\tilde X^{a}_{t})dt \Big],
	\end{equation*}
	with $\tilde X^{a}=X_{0}+\int_{0}^{\cdot} \sigma_{X}^{a_{t}}(\tilde X^{a}_{t})dW_{t}$. We therefore recover the right-hand side of \eqref{eq: dual intro} upon the change of variable $\alpha= \sigma_{X}^{a}(\tilde X^{a})$ $\Leftrightarrow$ $a=(\alpha-\sigma(X^{\alpha}))/f(X^{\alpha})$. 

	\vs2

In this paper, we consider a pretty general  path-dependent framework, generalizing the abstract setting of \cite{BLSZ18},  in which we  show that \eqref{eq: dual intro} always hold true whenever a solution to the dual optimal control problem exists (in a weak formulation approach) and we relate the solution of the dual problem to the one of the hedging problem. This shows that solving the second order FBSDE associated to the hedging problem boils down to solving a rather simple optimal control problem: its solution actually solves the highly non trivial  fixed point problem embedded in the 2FBSDE.  The link between the two solutions is simply given by the former change of variable  $a=(\alpha-\sigma(X^{\alpha}))/f(X^{\alpha})$, while the $b$ process is associated in an explicit manner to the (Fr\'echet) derivatives of the coefficients $\Phi$ and $(x,\alpha)\mapsto  |\alpha-\sigma(x)|^{2}/f(x)$ computed along the path of the optimal strategy. 
\vs2

In the Markovian setting of \cite{BLSZ18}, it is shown that the dual formulation \eqref{eq: dual intro} coincides with the minimum super-hedging price. 
	Therefore, $v(0,X_0)$ provides the cheapest hedging price. 
	On the other hand, by the above argument, the results in the current paper show that the dual formulation provides  the most expensive perfect hedging price/strategy.
	Combining the two implies that there is only one hedging price (and a unique hedging strategy, at least if $v$ is smooth enough), or, equivalently, only one solution to the 2FBSDE \eqref{eq: intro V Y}-\eqref{eq: intro X}.

\vs2 

	We shall not seek for a general existence result in the dual problem but  discuss only in Section \ref{sec: existence} some sufficient conditions, which are already far less restrictive than the ones imposed in \cite{BLSZ18}. 
	In fact, the existence of a solution to the dual optimal control problem \eqref{eq: dual intro} under weaker assumptions is itself non-standard.
	First, in general, the control problem is not concave. Second, the penalty term   in the dual formulation is only quadratic, while the usual weak existence results require a super-quadratic penalty to ensure $C$-tightness, see e.g.~\cite{haussmann1990existence}.
	We would like to leave it for future researches.

	\vspace{0.5em}
	
	{
	Finally, we   refer to \cite{Bilarev} for a variant market impact model with a resilience effect,
	and  to \cite{becherer2017stability, becherer2016optimal, becherer2016optimalbis, cetin2004liquidity, frey, Liu, Schon, Sircar}  for  related works on the hedging or liquidation  under market impact, see the introduction of \cite{BoLoZo1} for more details.
	}
	\vs2
	
	Before to conclude this introduction with notations and definitions that will be used all over this paper, let us also note that the course of our arguments leads us to formulate  a version of It\^{o}'s Lemma for  path-dependent functionals that are only $C^{0,1}$ in the sense of Dupire, which is of own interest, see the Appendix.

	\vspace{0.5em}

	{\bf Notations.} 
	$\mathrm{(i)}$. Let $\CT$ denote the space of all $\R$-valued continuous paths on $[0,T]$,
	and $\DT$ the space of all $\R$-valued c\`adl\`ag paths on $[0,T]$.	
	Given $(t,\xr)\in \TD$, we define the stopped path $\xr_{t \wedge \cdot}:=(\xr_{t\wedge s})_{s\in [0,T]}$,
	and   
	$$
		\Theta
		~:=~ 
		\big\{ (t,\xr)\in \TD ~: \xr=\xr_{t \wedge \cdot}  \big\}.
	$$
	For all $(t, \xr), (t', \xr') \in \TD$, let
	\begin{align*}
		\Ninfty{\xr- \xr'}
		&:=
		\sup_{s \in [0,T]}|\xr_s -  \xr'_s|
		\; \mbox{ and }\;
		\dr \big( (t,\xr), (t', \xr') \big)
		:=
		|t- t'|+ \Ninfty{\xr_{t \wedge \cdot} - \xr'_{t'\wedge \cdot}}.
	\end{align*}
	Notice that $\| \cdot \|$ defines a norm on $\DT$ and $\CT$, and that $\dr$ defines a distance on $\Theta$.
	
	\vspace{0.5em}

	$\mathrm{(ii)}$.  A function $\vp: \TD \to \R$ is said to be non-anticipative if $\vp(t, \xr)=\vp(t, \xr_{t \wedge \cdot})$ for all $(t,\xr)\in \TD$.
	It is clear that a non-anticipative function $f$ can be considered as a function defined on $\Theta$.
	We then denote by $\Cb(\Theta)$ the class of all non-anticipative functions $\vp: \Theta \longrightarrow \R$ such that
	$\vp(t^{n},\xr^{n})\to \vp(t,\xr)$ for all sequences $(t^{n},\xr^{n})_{n\ge 1}\subset \Theta$ and $(t,\xr)\in  \Theta$ satisfying $t^{n}\searrow t$ and 
	$\Ninfty{\xr^{n}-\xr}\to 0$.
	When $\DT$ is considered as the canonical space, a measurable function $\vp: \DT \to \R$ can be considered as a random variable,
	and a function $\vp: \TD \to \R$ can be considered as a process. We will therefore write $\vp_t(\xr)$ in place of $\vp(t, \xr)$.
	
	\vspace{0.5em}
	
	$\mathrm{(iii)}$.  Following \cite{cont2013functional}, a function $\vp: \Theta \to \R$ is said to be horizontally differentiable if,
	for all $(t,\xr)\in [0,T) \x \DT $, its horizontal derivative
	$$
		\partial_{t}\vp(t,\xr) ~:=~ \lim_{h\searrow 0} \frac{\vp(t+h,\xr_{t \wedge \cdot}) - \vp(t,\xr_{t \wedge \cdot})}{h}
	$$
	is well-defined.
	A function $\vp$ on $\Theta$ is said to be vertically differentiable if, for all $(t,\xr)\in \TD$, its vertical derivative
	$$
		\nabla_{\xr}\vp(t,\xr) ~:=~ \lim_{y\to  0,y\ne 0} \frac{\vp(t,\xr\oplus_{t} y)-\vp(t,\xr)}{y},
		~~\mbox{with}~~
		\xr\oplus_{t} y := \xr+\1_{[t,T]} y, ~\forall y \in \R,
	$$
	is well-defined.
	For a function $\vp: \Theta \longrightarrow \R$, 
	we say that $\vp\in \Cb^{0,1}(\Theta)$ if both $\vp$ and $\nabla_{\xr}\vp$ are well defined and belong to $\Cb(\Theta)$.
	We say that $\vp\in \Cb^{1,2}(\Theta)$ if $\vp$, $\nabla_{\xr}\vp$, $\nabla^{2}_{\xr}\vp:=\nabla_{\xr}(\nabla_{\xr}\vp)$  and $\partial_{t}\vp$ are all well defined and belong to $\Cb(\Theta)$.

	\vspace{0.5em}

	$\mathrm{(iv)}$. A function $\vp:\DT \longrightarrow \R$ is Fr\'echet differentiable if there exists 
	a measurable map $\xr \in \DT \longmapsto \lambda_{\vp}(\cdot,\xr)$ taking values in the space of finite measures on $[0,T]$, 
	and a measurable function $R_{\vp}: D([0,T])^{2} \longrightarrow \R$
 	such that, for each given $\xr\in \DT$, one has
	$$
		\vp( \xr') 
		~=~
		\vp(\xr)+\int_{0}^{T}( \xr'_{u}-\xr_{u})\lambda_{\vp}(du,\xr)+R_{\vp}(\xr,  \xr'), \;\xr'\in \DT,
	$$
	and $R_{\vp}( \xr, \xr')/\|\xr- \xr'\|\to 0$ as $\xr'\to \xr$. 
	When $\vp: \DT \longrightarrow \R$ is Fr\'echet differentiable, $\lambda_{\vp}$ is called the Fr\'echet derivative of $\vp$ and $R_{\vp}$ the residual term. 
	Finally, a non-anticipative map $\vp : \Theta \longrightarrow \R$ is said to be Dupire-concave if, 
	for all  $t \in [0, T]$ and $\xr^{1},\xr^{2}\in D([0,T])$ such that $\xr^{1}=\xr^{2}$ on $[0,t)$,
	one has
	$$
		\vp(t,\theta \xr^{1}+(1-\theta)\xr^{2})
		~\ge~ \theta \vp(t, \xr^{1})+(1-\theta)\vp(t,\xr^{2}),
		\;\;\mbox{for all}~\theta \in [0,1].
	$$
  
	$\mathrm{(v)}$. Denote by $\Om:= C([0,T])$ the canonical space of continuous paths on $[0,T]$,
	with canonical process $X$ and canonical filtration $\F$ generated by $X$.
	Notice that a process is $\F$-predictable if and only if it is $\F$-progressively measurable (see e.g. \cite[Proposition 9]{claisse2016pseudo}).
	Moreover, given a probability measure $\P$ on $\Om$, let $\F^{\P,+}$ be the $\P$-augmented filtration of $\F$.
	Then, a $\F^{\P,+}$-predictable process is $\P$-indistinguishable from a $\F$-predictable process (see e.g. \cite[Theorem IV. 78]{dellacherieMeyer}).
	We will then usually treat $\F^{\P,+}$-predictable processes as a predictable process w.r.t. the canonical filtration $\F$.
	
	\vspace{0.5em}
	
	$\mathrm{(vi)}$. Denote by $\Rb := \R \cup \{-\infty\}$,
	and by $\0$ the constant path in $\Om$ which equals   $0$ at every time $t \in [0,T]$.
	We use the convention that $\infty \times 0 = 0$, so that, for a constant process $(Z_r)_{r \in [s,t]}$, we have 
	$\int_s^t \infty dZ_r = 0$.
	 Finally, given  $0 \le s \le t \le T$, a function $f: [0,T] \to \R$ and a measure $\mu$ on $[0,T]$, 
	we write $\int_s^t f(r) \mu(dr) :=  \int_{[s,t]} f(r) \mu(dr)$ whenever the latter integration is well defined.

\section{Dual formulation and perfect hedging}

	In this section, we first define our primal hedging problem, which is an abstract version of the one considered in \cite{BLSZ18,BoLoZo2}. 
	We then define the corresponding dual problem and 
	{state our main result which provides an explicit construction of a solution to the primal hedging problem from a solution to the dual optimal control problem.}

\subsection{Primal hedging problem}

	We consider a weak formulation approach of the hedging problem presented in the introduction. The  dynamic of the underlying risky asset follows a diffusion process, but with path-dependent coefficients, and the payoff function is also path-dependent. Although the paths of the underlying asset are continuous, we shall need to  let the coefficient functions be defined on the space $D([0,T])$ of c\`adl\`ag paths, rather than on $\Om = C([0,T])$, to apply the It\^{o}'s functional calculus as well as calculus of variation arguments later on.

	\vs2
	
	Let $\Phi :  D([0,T])\mapsto \R$ be a measurable payoff function,
	$\sigma: [0,T]  \x D([0,T]) \x \Rb \to \R \cup \{ \infty\}$ and $F: [0,T] \x D([0,T]) \x \Rb \to \R \cup \{ \infty\}$ be measurable maps such that
	$(\sigma, F)(\cdot, \gamma)$ are non-anticipative for all $\gamma\in \Rb$.
	Let $x_0 \in \R$ be a fixed constant.
\vs2

	Given the above, our hedging problem is formulated as follows.
	\begin{definition} \label{def:hedge_pb}
		A probability measures $\P$ on $\Om$, together with 
		$$
			\mbox{two constants}~y_0, v_0 \in \R
			~\mbox{and}~
			\F \mbox{-predictable processes}~
			\big(Y, V, \Bk, \Gk \big),
		$$
		is a (weak) solution to the primal hedging problem if 

		\vspace{1mm}

		$\mathrm{(i)}$ $X$ is a $\P$-martingale such that $\P[X_0 = x_0] = 1$, $d \langle X \rangle_t = \sigma_t^2(X, \Gk_t) dt$ on $[0,T]$, $\P$-a.s., 
		and $\Gk = (\Gk_t)_{0 \le t \le T}$ is a $\Rb$-valued process satisfying
		$$
			\E^{\P} \Big[\int_{0}^{T}|\sigma_{t}(X, \Gk_{t})|^{2}dt \Big]<\infty
			~~\mbox{and}~~
			\int_{0}^{T}|\Gk_{t}|^{2}d\langle X \rangle_{t}<\infty,~~\P\mbox{-a.s.}
		$$
		
		$\mathrm{(ii)}$ $\Bk$ is {$\R$-valued with} bounded total variation, and $(Y, V)$ {are $\R$-valued and satisfy}:   
		$V_{T}=\Phi(X)$, $\P$-a.s, 
		\begin{equation} \label{eq:VY}
			V_t=v_0+ \!\int_{0}^{t} \!\! F_s(X, \Gk_s) ds +\! \int_{0}^{t} \!\! Y_s dX_s
			~\mbox{and}~
			Y_t = y_0 +\! \int_{0}^{t} \!\! \Gk_s dX_s -\Bk_t,
			~t \in [0,T],
			~\P\mbox{-a.s.}
		\end{equation}
	\end{definition}

	\begin{remark}
		Under the conditions of Definition \ref{def:hedge_pb}, it is well-known (see e.g.  \cite[Proposition 2.1, Chap.~IV]{ikeda2014stochastic})  that 
		there exists a possibly enlarged filtered probability space, equipped with a Brownian motion $W^{\P}$, such that
		\begin{equation} \label{eq:Dynamic_X}
			X_t = x_0 + \int_{0}^{t} \sigma_{s}(X, \Gk_{s}) dW^{\P}_{s},
			~~t \in [0,T],
			~\P\mbox{-a.s.}
		\end{equation}
		As mentioned in the introduction, \eqref{eq:VY}-\eqref{eq:Dynamic_X} is a {(non Markovian)} Second Order Backward Stochastic Differential Equation in the sens of \cite{cheridito2007second,soner2009dynamic}, with a specific  coupling    through the quadratic variation of $Y$ which appears both  in the quadratic variation of the forward process $X$, and in the drift of the backward process $V$. 
	\end{remark}

	\begin{remark}
		A solution $(\P, y_0, v_0, Y, V, B, \Gk)$ of the hedging problem   has the following financial interpretation.
		The probability $\P$ is a martingale probability measure,
		under which the underlying asset $X$ follows the dynamic \eqref{eq:Dynamic_X},
		with volatility function $\sigma$ impacted by $\Gk$. 
		The process $Y$ represents the dynamic trading strategy, i.e.~the number of units of the risky asset in the self-financing portfolio,
		whose value process is given by $V$ and which is also impacted by $\Gk$.
		The process $\Gk$ measures the sensitivity of the dynamic strategy $Y$ w.r.t. the evolution of the underlying asset $X$ (see \eqref{eq:VY}).
		The equality $V_{T}=\Phi(X)$ means that the path-dependent option with payoff function $\Phi$ is replicated perfectly by the portfolio $V$. 
	\end{remark}

	\begin{example}  \label{exam:BoLoZo}
		The model considered in  \cite{BoLoZo1,BoLoZo2},
		as described in the introduction, 
		corresponds to the following coefficients:
		\begin{eqnarray*}
			\sigma_t(\xr,\gamma)
			=
			\frac{ \sigma_0(t,\xr_t)}{1-f(\xr_t)\gamma} \1_{A} + \infty \1_{A^c},
			~~
			F_t(\xr,\gamma)
			=
			\frac12 \left(\frac{ \sigma_0(t,\xr_t){\gamma}}{1-f(\xr_t) \gamma}\right)^{2}f(\xr_t)
			\1_A +  \infty \1_{A^c},
		\end{eqnarray*}
		for some measurable functions $f: \R \to (0,\infty)$, $\sigma_0: [0,T] \x \R \to [0,\infty)$, and with 
		$A := \{(t,\xr,\gamma) ~: f(\xr_t) \gamma < 1 \}$. 
		In  \cite{BLSZ18,BoLoZo1,BoLoZo2},  the process $B$ is taken to be absolutely continuous. The dynamics \eqref{eq:VY} is a relaxation which aims at ensuring existence in a more general setting. Conditions under which $B$ will indeed be absolutely continuous are discussed in  Remark \ref{rem: B absolutly cont} below. Moreover, in  \cite{BoLoZo1,BoLoZo2}, conditions are imposed on the strategies, so that a (unique) martingale measure can be associated to each of them. As in the usual hedging literature, the drift of the price process does not play any role. In the current paper, we write the dynamics directly under the martingale measure. 
		This should also be viewed as a relaxation.
	\end{example}

\subsection{The dual problem}

	To define our dual problem, let us first introduce a technical condition on $\sigma$ and $F$,
	which is assumed throughout the paper.
	
	\begin{assumption} \label{assum:sigma_invertible}
		$\mathrm{(i)}$ For all $a \in (0, \infty)$ and every $(t,\xr)\in [0,T]\x   D([0,T])$, there exists a unique $\gamma \in \R$ such that $\sigma_t(\xr, \gamma) = a$, and we write  $\sigma^{-1}_t(\xr, a) := \gamma$.
		\vspace{0.5em}
		
		$\mathrm{(ii)}$ There exists a function $G: [0,T] \x \DT \x [0, \infty) \to \R$, such that, for all $(t,\xr) \in [0,T] \x \DT$ and $a > 0$,
		\begin{equation} \label{eq:structure}
			G_t(\xr, a) = F_t(\xr, \sigma^{-1}_t(\xr, a))
			~~\mbox{and}~~
			\partial_a G_t(\xr, a) = a \sigma^{-1}_t(\xr, a).
		\end{equation}
		Moreover, by setting $\sigma^{-1}_t(\xr, 0) := - \infty$, one has $ F_t(\xr, \sigma^{-1}_t(\xr, 0)) = F_t(\xr, -\infty) = G_t(\xr, 0)$.
	\end{assumption}

	\begin{remark} \label{rem:G}
		The conditions in Assumption \ref{assum:sigma_invertible} are mainly motivated by Example \ref{exam:BoLoZo}, 
		in which, by direct computation, one obtains that
		$$
			\sigma^{-1}_t(\xr, a)
			~=~
			\frac{a - \sigma_0(t, \xr_t)}{f(\xr_t) a},
			~~~
			G_t(\xr, a)
			~=~
			\frac12 \frac{(a-\sigma_0(t,\xr_t))^{2}}{f( \xr_t)},
			~~~
			\partial_a G_t(\xr, a) = \frac{a - \sigma_0(t,\xr_t)}{f(\xr_t)},
		$$
		for all $(t,\xr,a) \in [0,T] \x \DT \x [0, \infty)$.
	\end{remark}

	The dual problem is a standard  stochastic optimal control problem in a weak formulation (in the sense of e.g.~\cite{haussmann1990existence}), which can be formulated on the canonical space $\Om = C([0,T])$ of all continuous paths
	equipped with canonical process $X$ and canonical filtration $\F$.
	\vs2
	
	In order to define it properly, let us first  denote by $\Pc_0$ the collection of all probability measures  $\P$ on $\Om$, 
	under which $X$ is a square-integrable continuous martingale on $[0,T]$ such that $X_0 = 0$ a.s.~and its quadratic variation is absolutely continuous w.r.t. the Lebesgue measure a.s.
	Then, there exists (see e.g.  Karandikar \cite{karandikar1983quadratic}) a $\F$-predictable non-decreasing process $\langle X \rangle: [0,T] \x \Om \to \R \cup \{\infty\}$
	such that the process $\langle X \rangle$ is a version of the $\P$-quadratic variation of $X$ on $[0,T]$  for all $\P \in \Pc_0$, and one can define a $\F$-predictable process $ \alpha = ( \alpha_t)_{0 \le t \le T}$ by
	$$
		\alpha_t ~=~\alpha_{t}(X)~:=~
		\sqrt{ \limsup_{\varepsilon \searrow 0}
		~\frac{\langle X \rangle_t - \langle X \rangle_{(t-\varepsilon)\vee 0}}{\varepsilon}},
		~~\mbox{for all}~t \in [0,T].
	$$

		With this construction, one has
		$
			\langle X \rangle_t = \int_0^t \alpha^2_s ds,~t \in [0,T], ~\P\mbox{-a.s.},
		$
		for every $\P \in \Pc_0$.
		Moreover, on a possibly enlarged filtered probability space, there exists
		a Brownian motion $(W^{\P}_t)_{t \in [0,T]}$ such that
		\begin{equation} \label{eq:X_repres}
			X_t = X_0 + \int_0^t \alpha_s dW^{\P}_s, ~t \in [0,T],~\P\mbox{-a.s.}
		\end{equation}
		In other words, under $\P$, $X$ can be considered as a controlled diffusion process associated with the control process $\alpha$. 
		This leads to the following standard definitions of control rules.

	\begin{definition}\label{def : strook varadhan} 
		$\mathrm{(i)}$ A probability measure $\P \in \Pc_0$ is called a weak control rule (or simply a control rule).
		
		\vspace{0.5em}
		
		\noindent $\mathrm{(ii)}$  A weak control rule $\P \in \Pc_0$ is called a strong control rule 
		if there is a possibly enlarged probability space, together with the representation \eqref{eq:X_repres}, such that
		$(\alpha_t)_{t \in [0,T]}$ is equal to some process that is predictable w.r.t. the filtration generated by $(W^{\P}_t)_{t \in [0,T]}$,
		in the $d\P \x dt$-a.e. sens.
		We denote by $\Pc^S_0$ the collection of all strong control rules.
	\end{definition}

	Given $(t,\xr) \in \TD$, we now define the processes $\Xb^{t,\xr} = (\Xb^{t,\xr}_s)_{s \in [0,T]}$ and $\alphab^t = (\alphab^t_s)_{s \in [t,T]}$  on $\Om$   by
	$$
		\Xb^{t,\xr}_s(\om) ~:=~ \xr_{s \wedge t} + X_{(s-t) \vee 0}(\om) - X_0(\om)
		~~\mbox{and}~
		\alphab^t_s (\omega) ~:=~ \alpha_{(s-t) \vee 0}(\omega),
		~~\mbox{for all}~\om \in \Om.
	$$
	Observe that, under each $\P \in \Pc_0$,
	the process $\Xb^{t,\xr}$ has deterministic and c\`adl\`ag path $\xr$ as initial condition on $[0,t]$, and then evolves as a controlled continuous martingale diffusion process on $[t,T]$. 
	This construction allows us to restrict our optimization problem to the measures on the canonical space of continuous paths while allowing the path of the process before $t$ to be c\`adl\`ag (which will be useful from a pure technical point of view).
	Denote by $\Pcb_{t,\xr}$ (resp. $\Pcb^S_{t,\xr}$)  the collection of laws of $\Xb^{t,\xr}$ under the weak (resp. strong) control rules, i.e.
	$$
		\Pcb_{t,\xr} 
		~:=~ 
		\Big\{ \P \circ \big( \Xb^{t,\xr} \big)^{-1} 
			~:~
			\P \in \Pc_0
		\Big\}
		~~\mbox{and}~~
		\Pcb^S_{t,\xr}
		~:=~
		\Big\{ \P \circ \big( \Xb^{t,\xr} \big)^{-1} 
			~: 
			\P \in \Pc^S_0
		\Big\}.
	$$
	An element $\Pb \in \Pcb_{t,\xr}$ is a probability measure on $\DT$, as $\Xb^{t,\xr}$ has c\`adl\`ag paths in general.

	We can finally  introduce the value function of our dual problem:
	\begin{equation} \label{eq:dual_pb}
		\vr(t, \xr) ~:= \sup_{\P \in \Pc_0} J(t,\xr; \P),
		~~\mbox{where}~~
		J(t,\xr;\P) ~:=~ \E^{\P} \Big[\Phi\big( \Xb^{t,\xr} \big) - \int_{t}^{T}G_{s} \big( \Xb^{t,\xr}, \alphab^t_s \big) ds \Big].
	\end{equation}
	Denote also by $\Pcb^*_{t,\xr}$ the set of laws of the controlled process under the optimal (weak) control rules, i.e.
	$$
		\Pcb^*_{t,\xr} 
		~:=~ 
		\Big\{
			\P \circ \big( \Xb^{t,\xr} \big)^{-1} 
			~:
			\P \in \Pc_{0},
			~J(t,\xr; \P) = \vr(t, \xr)
		\Big\}.
	$$

 	\begin{remark} \label{rem:Pc}
		$\mathrm{(i)}$
		If $\xr_{t \wedge \cdot} \in \Om = C([0,T])$, so that $\Xb^{t,\xr}$ is a process with continuous paths,
		then $\Pcb_{t,\xr}$, $\Pcb^S_{t,\xr}$ and $\Pcb^*_{t,\xr}$ can also be considered as sets of probability measures on $\Om$.
		In this case, we simply write
		$$
			\Pc_{t,\xr}, ~\Pc^S_{t,\xr} ~\mbox{and}~\Pc^*_{t,\xr},\mbox{ as sets of probability measures on $\Om$, in place of}~\Pcb_{t,\xr}, ~\Pcb^S_{t,\xr}~\mbox{and}~\Pcb^*_{t,\xr}.
		$$
		In particular,  when $\xr_{t \wedge \cdot} \in C([0,T])$, one has the equivalent formulation of the control problem:
		\begin{equation} \label{eq:dual_pb_equiv}
			\vr(t, \xr) = \sup_{\P \in \Pc_{t,\xr}} \Jb(t, \xr; \P),
			~~\mbox{where}~~
			\Jb(t,\xr;\P) ~:=~ \E^{\P} \Big[\Phi\big( X \big) - \int_{t}^{T}G_{s} \big( X, \alpha_s \big) ds \Big],
		\end{equation}
		and $\Pc^*_{t,\xr}$ turns out to be the set of the all optimizers for the optimization problem \eqref{eq:dual_pb_equiv}.

		\vspace{0.5em}

		$\mathrm{(ii)}$ 
		Notice also that $\vr(t, \xr)$ (resp. $\Pcb_{t,\xr}$, $\Pcb^*_{t,\xr}$) is non-anticipative in the sense that $\vr(t,\xr) = \vr(t, \xr_{t \wedge \cdot})$ (resp. $\Pcb_{t,\xr} = \Pcb_{t,\xr_{t\wedge \cdot}}$, $\Pcb^*_{t,\xr} = \Pcb^*_{t,\xr_{t\wedge \cdot}}$).
		In particular, when $t=0$, these functions and sets depend only on $x_0 := \xr_0 \in \R$.
		We then simplify their notations to $\vr(0, x_0)$, $\Pcb_{0,x_0}$, $\Pcb^*_{0,x_0}$, $\Pc_{0, x_0}$, $\Pc^*_{0, x_0}$, $\Pc^S_{0,x_0}$, etc.
	\end{remark}

 \subsection{Main result}
 
	Our main result consists in providing an explicite construction of a solution to the  primal hedging problem in Definition \ref{def:hedge_pb}, based on an optimal solution to the dual problem \eqref{eq:dual_pb}.
	\vs2
	
	Besides Assumption \ref{assum:sigma_invertible}, we shall assume some additional technical conditions on the functionals $\Phi: \DT \to \R$ 
	and $G: [0,T] \x \DT \x [0,\infty) \to \R$.

	\begin{enumerate}[\rm (H1).]
		\item  \label{ass:Phi_Frechet_diff}
		The function $\Phi$ is Fr\'echet differentiable with Fr\'echet derivative  $\lambda_{\Phi}(\cdot,\xr)$,
		which is bounded for the total variation norm, uniformly in $\xr\in D([0,T])$.

		\item\label{ass:G_Frechet_diff} 
			The function $G$ is continuous on $[0,T] \x \DT \x [0,\infty)$.
			For every $a \in \R_+$, the function $(t,\xr) \mapsto G_t(\xr, a)$ is non-anticipative.
			Further, for every $(t,a)\in [0,T]\x \R$, $G_{t}(\cdot,a)$ is Fr\'echet differentiable with Fr\'echet derivative  $\lambda_{G}(\cdot,\xr;t,a)$. 
			The map $(t,\xr,a)\in \TD\x \R \mapsto \lambda_{G}(\cdot,\xr;t,a)$ is measurable. 
			Moreover, there exists a measurable map $\partial_{a}G :  \TD\x \R\to \R$ and $C>0$ such that, for all given $(t,\xr,a) \in [0,T]\x \DT\x \R$,  
			\begin{align*}
				&G_{t}(\bar \xr,\bar a)~=~G_{t}( \xr, a)+\int_{0}^{T}(\bar \xr_{u}-\xr_{u})\lambda_{G}(du,\xr;t,a)+\partial_{a}G_{t}(\xr,a)(\bar a-a) \\
				&~~~~~~~~~~~~~~~~~~~~~~~~~~~~~~~+o(\|\bar\xr-\xr\|+|\bar a-a|),\;\;\;\;\;\mbox{ 	for all   $(\bar \xr,\bar a)\in \DT\x \R$,}\\
				&\mbox{and } \;
				|\partial_{a}G_{t}(\xr,a)|+|\lambda_{G}|([0,t],\xr;t,a)\le C(1+|a|^{2}).
			\end{align*}

		\item\label{ass:bound_G} 
		For some constants $C>0$ and $C_0 > 0$, one has
		$$
			\frac{a^2}{C} - C
			~\le~
			G_t(\xr, a)
			~\le~
			C_0(1+a^2),\;\mbox{for all $(t,\xr,a)\in \TD\x \R_+$.}
		$$

		\item\label{ass:concavity} \label{ass:last}
		For all $t < T$ and $M > 0$, one has  
		$$
			\lim_{h \searrow 0} 
			\Big( \sup_{\bar \xr \in D^M_{t+h}} 
			\Big|
				\int_t^{t+h}  \!\!\! \lambda_{\Phi} (du,  \bar \xr)
			\Big|	
		+ 
		\sup_{\P\in \widetilde \Pc_0^M} 
		\E^{\P}\Big[
		\sup_{\bar \xr \in D^M_{t+h}}
			\Big|
				\int_{t+h}^T \int_t^{t+h} \!\!\! \lambda_G(du, \bar \xr; s, \alphab^{t+h}_s) ds
			\Big|
		\Big]
		\Big)
		~=~0,
		$$
		where
		\begin{equation} \label{eq:def_DM}
			D^M_{t+h} ~:=~ \big\{
			\bar \xr \in \DT
			~: \| \bar \xr_{t +h \wedge \cdot} \| \le M
			\big\}
			~\mbox{and}~
			\Pct_0^M := \Big\{ \P \in \Pc_0 ~:  \E^{\P} \Big[ \int_0^T\alpha^2_s ds \Big] \le M  \Big\}.
		\end{equation}

	\end{enumerate}

 	\begin{example}
		Let $\mu$ be a probability measure on $[0,T]$ without atom. Given $\phi \in \C^1_b(\R^2)$,
		   let $\Phi$ be defined by
		$$
			\Phi(\xr) ~=~ \phi \big( \mathrm{a}(\xr),  \xr_T \big),
			~~\mbox{with}~~
			 \mathrm{a}(\xr) := \int_0^T \xr_t \mu(dt).
		$$
		It follows from direct computations that
		$$
			\lambda_{\Phi}(dt, \xr) 
			~=~
			\partial_a \phi( \mathrm{a}(\xr), \xr_T)  \mu(dt)
			~+~
			\partial_x \phi( \mathrm{a}(\xr), \xr_T) \delta_T(dt).
		$$
		Let now $G$ be given as in  Remark \ref{rem:G}, for some bounded and positive smooth functions $\sigma_0$ and $f$.
		Then, for all $\xr \in \DT$, $\lambda_G([t, t+h], \xr; s, a) = 0$ whenever $s \in (t+h, T]$, so that
		$$
			\int_{t+h}^T \int_t^{t+h} \!\!\! \lambda_G(du, \bar \xr; s, \alphab^{t+h}_s) ds = 0.
		$$
		Based on the above computation, it is easy to check that Conditions \Href{ass:Phi_Frechet_diff}-\Href{ass:last} hold true.	
	\end{example}
 
	\begin{remark}\label{rem: lin growth when bounded frechet derivative}  
		Let a function $\Xi$ be Fr\'echet differentiable with Fr\'echet derivative  $\lambda_{\Xi}(\cdot,\xr)$,
		which is bounded for the total variation norm, uniformly in $\xr\in D([0,T])$. Then, $\Xi$ has at most linear growth: there exists $C>0$ such that 
		$$
		|\Xi(\xr)|\le C(1+\|\xr\|), \;\mbox{ for all }\xr \in \DT.
		$$ 
		Indeed, as $\partial_{r}\Xi(r\xr)=\int_{0}^{T}  \xr_{t} \lambda_{\Xi}(dt,  r\xr)$, for $r\in [0,1]$,
		one obtains the above inequality from the uniform boundedness of the total variation of  $\lambda_{\Xi}$. Similarly,  
		$$
		\Xi(\xr')=\Xi(\xr)+\int_{0}^{1}\int_{0}^{T} (\xr'_{t}-\xr_{t}) \lambda_{\Xi}(dt, r(\xr'-\xr)+\xr) dr,\; \mbox{  $\xr,\xr'\in \DT$.}
		$$
	\end{remark}
	
	Let us now introduce  two important processes that will play the role of optimal controls in our primal problem. 
	First, let us consider the measurable process $A = (A_t)_{0 \le t \le T}$ on $\Om$ defined by
	\begin{equation} \label{eq:def_At}
		 A_t :=
		\int_{[0,t]} \lambda_{\Phi}(du, X)-\int_{[0,t]}   {\widehat \lambda_G}(du ),
		~t \in [0,T]
	\end{equation}
 where $ \widehat \lambda_G  $ is the measure defined by
		$$
			  \widehat \lambda_G (I ) := \int_0^T \int_0^T \1_{[0,s] \cap I}(u) \lambda_G(du, X; s, \alpha_s) ds, ~~\mbox{for all Borel sets}~I \subseteq [0,T].
		$$
 	For later use, note that 
 	\begin{align}\label{eq: AT -At}
	A_{T}-A_{t-}=\int_{t}^{T} \lambda_{\Phi}(du, X)-\int_t^T \int_t^s \lambda_G(du, X; s, \alpha_s) ds, \;t\le T,
	\end{align}
	with $A_{0-}=0$.

	\begin{remark}\label{rem:dual_proj}
		 The process $A = (A_t)_{0 \le t \le T}$ is  of $\P$-integrable variation, for all $\P \in \Pc_0$.
		Indeed, by Conditions \Href{ass:Phi_Frechet_diff} and \Href{ass:G_Frechet_diff}, there is a constant $C$ such that, for all $(s, \xr, a) \in [0,T] \x \DT \x [0, \infty)$,
		$$
			\big|\lambda_{\Phi} \big| ([0,T], \xr) \le C
			~~\mbox{and}~~
			\big| \lambda_{G} \big|([0,T], \xr; s, a) \le C(1 + a^{2}).
		$$
		Then, by direct computations, one has
		$$
			d\big| A_t \big| 
			~\ll~
			\big|\lambda_{\Phi} \big| (dt, X)
			~+~
			2 \big| \lambda_{G}  ([0,t],X; t,  {\alpha_t}) \big| dt
			~+~
			 {\big| \widehat \lambda_G \big| (dt)},
		$$
		 {where} $\alpha$ is $\P$-a.s.~square integrable under each $\P\in \Pc_{0}$.
		 Consequently, for every $\P \in \Pc_0$, 
		the process $A$ has a $(\P, \F)$-dual predictable projection, hereafter 
		denoted by $B^{\P} = (B^{\P}_{t})_{0 \le t \le T}$ (see \cite[Definition VI.73]{dellacherie1982probabilities} for a precise definition of the dual predictable projection).
	\end{remark}
	Let us finally define the $\F$-predictable process $\gammah = (\gammah_t)_{t \in [0,T]}$ by 
	\begin{equation} \label{eq:def_gamma}
		\gammah_t 
		~:=~
		\sigma^{-1}_t(X, \alpha_t)
		~=~
		\frac{\partial_a G_t(X, \alpha_t)}{\alpha_t} \1_{\{\alpha_t > 0\}} 
		~-~
		\infty \1_{\{\alpha_t = 0 \}}.
	\end{equation}
	
	Now we are in a position to state our main result.

	\begin{theorem}\label{thm:dualty} 
		Let Assumption \ref{assum:sigma_invertible} and \Href{ass:Phi_Frechet_diff}-\Href{ass:last}  hold true.
		Assume in addition that $\Pcb^*_{t,\xr}$ is nonempty for all $(t,\xr)\in \DT$. Then,\\
		$\mathrm{(i)}$ 	The value function $\vr \in  \Cb^{0,1}(\Theta)$.\\
		$\mathrm{(ii)}$ For each $x_0 \in \R$, there exists an optimal solution $\Ph \in \Pc^*_{0,x_0}$ to the dual control problem \eqref{eq:dual_pb_equiv} with initial condition $(0,x_0)$,
		such that the measure $\Ph$, together with the constants $(\vr(0,x_0), \nabla_{\xr} \vr(0,x_0))$ 
		and the processes
		$$
			\big(\nabla_{\xr} \vr(t, X_{t\wedge \cdot}), \vr(t, X_{t\wedge \cdot}), \Bk^{\Ph}_{{t-}}, \gammah_{t} \big)_{t \in [0,T]},
		$$
		provides a weak solution to the hedging problem of Definition \ref{def:hedge_pb}.	
 	\end{theorem}
	
	We postpone to Section \ref{sec: existence} the discussion of sufficient conditions for $\Pcb^*_{t,\xr}$ to be nonempty for all $(t,\xr)\in \DT$.
	
	\begin{remark}
		Notice that the process $\gammah$ defined by \eqref{eq:def_gamma} is generally $\R \cup\{-\infty\}$-valued.
		Nevertheless, for every $\P \in \Pc_0$, one has $\alpha > 0$, $d\P \x d\langle X \rangle$-a.e.,
		which implies that $\gammah$ is in fact $\R$-valued $d\P \x d\langle X \rangle$-a.e.
		When $\alpha  > 0$, $d \P \x dt$-a.e.,
		the process $\gammah$ is $\R$-valued under $\P$.
	\end{remark}
  
	\begin{remark}\label{rem: B absolutly cont} 
	 Note that $\Bk^{\Ph}$ can be ensured to be absolutely continuous on the half-open interval $[0,T)$  by imposing suitable assumptions on $\lambda_{\Phi}$  and $\lambda_{G}$, to match with the setting of \cite{BoLoZo1,BoLoZo2}.
	\end{remark}

\section{Study of the dual optimal control problem}

 	The proof of Theorem \ref{thm:dualty} relies on a suitable characterization of an optimal control rule for \eqref{eq:dual_pb},
	together with a verification argument, which requires to apply the functional It\^o calculus of Proposition \ref{prop: Dupire Tanaka corrige} in the Appendix on the value function $\vr$, see Proposition \ref{prop : ito surbar v}  and Section \ref{sec: construction of the hedging strat} below.
	This requires to establish first that $\vr$ satisfies certain   regularity and concavity properties, and that  our dual optimal control problem is time consistent.

	\vspace{0.5em}
	
 	All over this section, we assume the conditions of Theorem \ref{thm:dualty}. In particular,   Assumption \ref{assum:sigma_invertible} and \Href{ass:Phi_Frechet_diff}-\Href{ass:last}  are in force. From now on, $C>0$ denotes a generic constant that does not depend on $(t,\xr)$, unless something else is specified.

\subsection{A-priori estimates, dynamic programming principle and time consistency of an optimal control}

	Let us first provide some basic properties of the reward function $J$ and the value function $\vr$ that will be used in the proof of the dynamic programming principle of Proposition \ref{prop:equiv_dpp}.
	
	\begin{proposition} \label{prop:basic_prop}
		For every $\P \in \Pc_0$, the reward function $J(\cdot; \P)$ is non-anticipative and belongs to $\Cb(\Theta)$.
		The value function $\vr$ is non-anticipative and there exists a  constant $C_{\vr}>0$ such that 
		\begin{equation} \label{eq:growth_v}
			|\vr(t, \xr)| 
			~\le~
			C_{\vr}(1 + \|\xr_{t\wedge \cdot}\|),
			~~\mbox{for all}~
			(t,\xr) \in \Theta.
		\end{equation}
	\end{proposition}
	\proof 
	Let $\P \in \Pc_0$ be fixed.
	As $G(\cdot, a)$ is non-anticipative, it follows directly that $J(\cdot; \P)$ and $\vr$ are both non-anticipative.
	
	\vspace{0.5em}
	
	Further, let $(t_n, \xr_n) \in \Theta$ be such that $t_n \ge t$ and $(t_n, \xr_n) \to (t,\xr)$ 
	(i.e. $ \| \xr_n(t_n \wedge \cdot) - \xr(t\wedge \cdot) \| \to 0$).
	Set $G_s(\cdot) \equiv 0$ for $s > T$,  and notice that $\alphab^{t}_s = \alphab^{t_n}_{s-t+t_n} = \alpha_{s-t}$ for $s \ge t$,
	so that
	$$
		\int_{t_n}^T G_s(\Xb^{t_n, \xr_n}, \alphab^{t_n}_s) ds 
		~=~ 
		\int_t^T G_{t_n -t+s}(\Xb^{t_n, \xr_n}, \alphab^t_s) ds. 
	$$
	Since $\| \Xb^{t_n, \xr_n}_{\cdot} - \Xb^{t, \xr}_{\cdot} \| \to 0$,
	and since $\Phi$ and $G$ are   continuous, it follows that
	$$
		\Phi(\Xb^{t_n, \xr_n}) - \int_{t_n}^T G_s(\Xb^{t_n, \xr_n}, \alphab^{t_n}_s) ds 
		~~\longrightarrow~~
		\Phi(\Xb^{t, \xr}) - \int_{t}^T G_s(\Xb^{t, \xr}, \alphab^{t}_s) ds,
		~~~\mbox{as}~~
		n \longrightarrow \infty.
	$$
	Because there exists $C>0$ such that $|\Phi(\xr)| \le C(1 + \|\xr\|)$ and $|G_s(\xr, a)| \le C(1 + a^2)$, for all $(s,\xr,a)\in [0,T]\x D([0,T])\x \R$, by \Href{ass:Phi_Frechet_diff}, Remark \ref{rem: lin growth when bounded frechet derivative} and \Href{ass:bound_G}, it is thus enough to use the dominated convergence theorem to deduce that $J(t_n, \xr_n; \P) \to J(t, \xr; \P)$,
	and therefore that $J(\cdot; \P) \in \Cb(\Theta)$.
	
	\vspace{0.5em}

	We finally  prove \eqref{eq:growth_v}. 
	By considering the weak control rule $\P \in \Pc_0$ satisfying $\P[X = \0] = 1$, one deduces from \Href{ass:bound_G}   that
	$$
		\vr(t, \xr) ~\ge~ J(t, \xr; \P) ~\ge~ \Phi(\xr_{t \wedge \cdot}) - C T~\ge~ -C' (1+T + \| \xr_{t \wedge \cdot}\|),
	$$
	for some $C'>0$ independent of $(t,\xr)$.
	Using again \Href{ass:Phi_Frechet_diff}, Remark \ref{rem: lin growth when bounded frechet derivative} and \Href{ass:bound_G}, and the fact that $\E^{\P}[ \int_0^t \alpha^2_s ds] = \E^{\P}[ X_t ^2]$ for all $\P \in \Pc_0$, 
	it follows that
	$$
		\vr(t, \xr) \le \sup_{\P \in \Pc_0} \Big( C\|\xr_{t \wedge \cdot}\| + C + C \E^{\P} [\|X_{(T-t) \wedge \cdot}\|] + C(T-t) - \frac{1}{C} \E^{\P} \big[ X_{T-t}^2 \big]  \Big)
		\le
		C_{\vr} (1+ \| \xr_{t \wedge \cdot} \|),
	$$
	for some $C_{\vr} >0$ independent of $(t,\xr)$.
	\endproof
	
	\begin{remark} \label{rem:optim_ctrl_bound} 
		For later use, note that it follows from  \eqref{eq:growth_v} that,
		in the optimal control problem \eqref{eq:dual_pb} with initial condition $(t, \xr)$,
		it is enough to consider $\P \in \Pc_0$ such that
		$$
			\E^{\P} \Big[ \Phi(\Xb^{t,\xr}) - \int_t^T G_s(\Xb^{t,\xr}, \alphab^t_s) ds \Big]
			~\ge~
			- C_{\vr} (1 + \|\xr_{t \wedge \cdot} \|).
		$$
		Moreover, by Conditions \Href{ass:Phi_Frechet_diff}, Remark \ref{rem: lin growth when bounded frechet derivative} and \Href{ass:bound_G}, the above inequality implies that
		\begin{equation*}
			\E^{\P} \Big[ \int_t^T\alphab^t_s ds \Big] ~\le~  C \sqrt{  (1+\|\xr_{t \wedge \cdot}\|)},
		\end{equation*}
		for some $C>0$ independent on $(t,\xr)$. 
		In other words, one has
		$$
			\vr(t, \xr) ~=~ \sup_{\P \in \Pct_0^M} J(t, \xr; \P), \;\forall\;(t,\xr)\in [0,T]\x\DT,
		$$
		for $M>0$ large enough,  
		where $\Pct_0^M$ is defined in \eqref{eq:def_DM}.
 	\end{remark}

	\vspace{0.5em}

	We next provide a classical dynamic programming result for the control problem \eqref{eq:dual_pb}.

	\begin{proposition} \label{prop:equiv_dpp}
		The value function $\vr$ is universally measurable and satisfies
		\begin{equation} \label{eq:DPP}
			\vr(t,\xr)
			~=~
			\sup_{\P \in \Pc_0}
			\E^{\P} \Big[
				\vr(\tau, \Xb^{t,\xr}_{\tau \wedge \cdot}) - \int_{t}^{\tau} G_{s}(\Xb^{t,\xr}, \alphab^t_{s}) ds
			\Big],
		\end{equation}
		for all $(t,\xr) \in \TD$ and all $\F$-stopping times $\tau$ taking value in $[t,T]$.
	\end{proposition}
	\proof For   $n \ge 1$, let us define
	$$
		\vr^n(t, \xr)
		:=
		\sup_{\P \in \Pc^n_0} J(t,\xr; \P)
		~\mbox{with}~\Pc^n_0 := \big\{ \P \in \Pc_0 ~:  \alpha_s \le n, ~d\P \x ds \mbox{-a.e.} \big\}.
	$$
	Then, $\vr^n$ is universally measurable, and it satisfies the dynamic programming principle (see e.g. \cite[Corollary 3.7]{el2013capacities}):
	\begin{equation} \label{eq:DPP_n}
		\vr^n(t,\xr)
		~=~
		\sup_{\P \in \Pc^n_0}
		\E^{\P} \Big[
			\vr^n(\tau, \Xb^{t,\xr}_{\tau \wedge \cdot}) - \int_{t}^{\tau} G_{s}(\Xb^{t,\xr}, \alphab^t_{s}) ds
		\Big].
	\end{equation}
	Moreover, for every $\P \in \Pc_0$ and $(t,\xr) \in \TD$, one can appeal to the representation \eqref{eq:X_repres} 
	and define
	$$
		\Xb^{t,\xr,n}_s := \xr_{s\wedge t}+ \int_0^{(s-t) \vee 0} \alpha^n_r dW^{\P}_r,
		~~\mbox{with}~~
		\alpha^n_r := \alpha_r \wedge n.
	$$
	Then,
	$$
		\E^{\P} \Big[ \sup_{s \in [0,T]} \big| \Xb^{t,\xr}_s - \Xb^{t,\xr,n}_s \big|^2 \Big]
		+
		\E^{\P} \Big[ \int_0^T \big( \alpha_s - \alpha^n_s \big)^2 ds \Big]
		~\longrightarrow~
		0.
	$$
	Since $\Phi$ has at most of linear growth, by \Href{ass:Phi_Frechet_diff} and Remark \ref{rem: lin growth when bounded frechet derivative}, and since
	$$
		\E^{\P} \Big[ \int_t^T \Big| G_s(\Xb^{t, \xr}, \alpha_s) - G_s(\Xb^{t,\xr, n}, \alpha^n_s) \Big| \1_{\{ \alpha_s > R\}} ds \Big]
		~\le~
		2 \E^{\P} \Big[ \int_t^T C(1+ \alpha_s^2) \1_{\{ \alpha_s > R\}} ds \Big]
		~\longrightarrow~
		0,
	$$
	as $R \to \infty$, recall \Href{ass:G_Frechet_diff} and  the definition of $\alpha^{n}$ above, it is   easy to deduce that
	$$
		J(t, \xr; \P^n) \longrightarrow J(t,\xr; \P),
		~\mbox{as}~ n \longrightarrow \infty,
		~\mbox{with}~
		\P^n := \P \circ \Big( \int_0^{\cdot} \alpha^n_s dW^{\P}_s \Big)^{-1} \in \Pc^n_0.
	$$
	This implies that $\vr^n(t, \xr) \nearrow \vr(t,\xr)$ as $n \longrightarrow \infty$. Hence,
	 $\vr$ is also universally measurable.
	
	\vspace{0.5em}
	
	Next, for $n \ge m \ge 1$,
	combining the dynamic programming principle \eqref{eq:DPP_n}  and the fact that 
	$\Pc^m_0 \subset \Pc^n_0 \subset \Pc_0$ 
	leads to
	$$
		\sup_{\P \in \Pc_0^m}
		\E^{\P} \Big[
			\vr^n(\tau, \Xb^{t,\xr}_{\tau \wedge \cdot}) - \int_{t}^{\tau} G_{s}(\Xb^{t,\xr}, \alphab^t_{s}) ds
		\Big]
		\le
		\vr^n(t,\xr)
		\le
		\sup_{\P \in \Pc_0}
		\E^{\P} \Big[
			\vr(\tau, \Xb^{t,\xr}_{\tau \wedge \cdot}) - \int_{t}^{\tau} G_{s}(\Xb^{t,\xr}, \alphab^t_{s}) ds
		\Big].
	$$
	Taking the limit $n \to \infty$ and using the monotone convergence theorem, we obtain
	$$
		\sup_{\P \in \Pc_0^m}
		\E^{\P} \Big[
			\vr(\tau, \Xb^{t,\xr}_{\tau \wedge \cdot}) - \int_{t}^{\tau} G_{s}(\Xb^{t,\xr}, \alphab^t_{s}) ds
		\Big]
		\le
		\vr(t,\xr)
		\le
		\sup_{\P \in \Pc_0}
		\E^{\P} \Big[
			\vr(\tau, \Xb^{t,\xr}_{\tau \wedge \cdot}) - \int_{t}^{\tau} G_{s}(\Xb^{t,\xr}, \alphab^t_{s}) ds
		\Big].
	$$
	Finally, as $|\vr(t, \xr)| \le C_{\vr}(1+ \|\xr\|)$ by Proposition \ref{prop:basic_prop},
	it is enough to use again the same arguments as in the proof of the assertion $\vr^n(t, \xr) \nearrow \vr(t,\xr)$ and then send $m \to \infty$ to deduce \eqref{eq:DPP}.
	\endproof

	\begin{remark}
		At this stage, we simply assert that the value function $\vr$ is universally measurable in Proposition \ref{prop:equiv_dpp}.
		In fact, it  will be proved to have much more regularity later in Section \ref{subsec:continuity}.
	\end{remark}

	\begin{lemma}\label{lemm:optimal_ctrl_cond}  
		Fix $x_0 \in \R$ and let $\Ph \in \Pc_{0,x_0}$ be an optimal control rule for \eqref{eq:dual_pb} with initial condition $(0, x_0)$.
		Given a stopping time $\tau$ taking values in $[0,T]$, one can consider a family $(\Ph_{\om})_{\om \in \Om}$ of regular conditional probability distributions  of $\Ph$ knowing $\Fc_{\tau}$.
		Then, for $\Ph$-a.e. $\om \in \Om$, the probability $\Ph_{\om}$ is an optimal solution to control problem \eqref{eq:dual_pb_equiv} with initial condition $(\tau(\om), X_{\tau(\omega)\wedge \cdot}(\om))$.
	\end{lemma} 
	\proof This is a direct consequence of the dynamic programming principle of Proposition \ref{prop:equiv_dpp}.
	Indeed,  it implies that
	\begin{eqnarray*}
		\E^{\Ph}\Big[  \vr(\tau,  X) - \int_{0}^{\tau}G_{s}( X,\alpha_{s})ds\Big]
		&\le&
		\vr(0, x_0)
		~=~
		\E^{\Ph} \Big[ \E^{\Ph} \Big[ \Phi(X) - \int_{0}^T G_{s}( X, \alpha_{s})ds \Big| \Fc_{\tau} \Big]  \Big]\\
		&=&
		\E^{\Ph} \Big[  Z - \int_{0}^{\tau}G_{s}( X, \alpha_{s})ds \Big],
 	\end{eqnarray*}
	where $Z: \Om \to \R$ is defined by
	$$
		Z(\om) ~:=~ \overline J(\tau(\om), X(\om); \Ph_{\om}).
	$$
 	On the other hand,  for $\Ph$-a.e. $\om \in \Om$, one has $Z(\om) = \overline J(\tau(\om), X(\om);  \Ph_{\om} ) \le  \vr(\tau(\om),  X(\om))$ by \eqref{eq:dual_pb_equiv}. 
	It follows that $Z(\om) = \overline J(\tau(\om), X(\om); \Ph_{\om}) = \vr(\tau(\om) X(\om))$, for $\Ph$-a.e. $\om \in \Om$.
	\endproof

\subsection{Continuity and concavity of $\vr$}
\label{subsec:continuity}

	We now establish two important technical results. 
	The first one is on the Dupire horizontal derivative of the gain function $J$ in \eqref{eq:DJ}. 
	Although it is elementary, combined with the {\sl enveloppe theorem} (Proposition \ref{prop : regu et decroissance de bar v}) below, it will play a major role in the proof of Proposition \ref{prop: characterization hat s}. 
	Next, we derive a time monotonicity and a concavity property for the value function $\vr$, 
	that will allow us to apply the functional It\^{o}'s formula of Proposition \ref{prop: Dupire Tanaka corrige} in the Appendix, but also to prove the {\sl enveloppe theorem} mentioned above. 

	\begin{proposition}\label{prop : continuity of J} 
		$\mathrm{(i)}$ For every $\P \in \Pc_0$, one has
		\begin{equation} \label{eq:DJ}
			\nabla_{\xr} J(t,\xr; \P)
			~=~
			\E^{\P} \Big[ 
				\int_t^T \lambda_{\Phi}(ds,\Xb^{t,\xr})
				+
				\int_t^T \int_t^s \lambda_G(du, \Xb^{t,\xr}; s, \alphab^t_s) ds
			\Big].
		\end{equation}
		
		$\mathrm{(ii)}$ 
		The value function $\vr \in \Cb(\Theta)$,
		and there exists a constant $C> 0$ such that
		\begin{equation} \label{eq:non_increas} 
			\vr(t+h, \xr_{t \wedge \cdot}) -  \vr(t, \xr)  ~\le~ Ch,
		\end{equation}
		for all $(t,\xr) \in \Theta,~h \in [0,T-t]$.
	\end{proposition}
	\proof $\mathrm{(i)}$ For $\P \in \Pc_0$ fixed, 
	the vertical derivative $\nabla_{\xr}J(\cdot;\P)$ can be computed directly under Conditions \Href{ass:Phi_Frechet_diff} and \Href{ass:G_Frechet_diff}  (also  recall Remark \ref{rem: lin growth when bounded frechet derivative} which allows to use the dominated convergence theorem under the above assumptions).
	
	\vspace{0.5em}
	
	$\mathrm{(ii)}$ From Remark \ref{rem:optim_ctrl_bound}, we know that
	$$
		 \vr(t,\xr) = \sup_{\P \in \widetilde \Pc_0^M}  J(t, \xr; \P),
		 ~~\mbox{for}~M~\mbox{large enough},
	$$
	where $\widetilde \Pc_0^M$ is defined in \eqref{eq:def_DM}.
	Let us fix $t \in [0,T)$, $\xr, \xr' \in \DT$ such that $\xr_{t\wedge \cdot} \neq \xr'_{t\wedge \cdot}$. Then, it
 follows from \eqref{eq:DJ}, combined with Conditions \Href{ass:Phi_Frechet_diff} and \Href{ass:G_Frechet_diff}, that
	\begin{equation} \label{eq:v_Lip}
		\frac{|\vr(t,\xr) - \vr(t, \xr') |}{ \|\xr_{t\wedge \cdot} - \xr'_{t \wedge \cdot} \|}
		~\le~
		C \Big(1+ \sup_{\P \in \widetilde\Pc_0^M} \E^{\P} \Big[  \int_t^T |\alphab^t_s|^2 ds \Big] \Big)
		~<~
		\infty,
	\end{equation}
	for some $C>0$ (independent of $(t,\xr,\xr')$).
	
	\vspace{0.5em}
	
	Next, recall  that $J(\cdot; \P) \in \Cb(\Theta)$ for all $\P \in \Pc_0$, by Proposition \ref{prop:basic_prop}, so that  there exists an optimal control $\Ph \in \Pc_0$ such that
	$$
		\vr(t, \xr) 
		~=~
		J(t, \xr; \Ph)
		~=~
		\lim_{t' \searrow t} J(t', \xr_{t \wedge \cdot}; \Ph) 
		~\le~
		\liminf_{t' \searrow t} \vr(t', \xr_{t \wedge \cdot}).
	$$
	
	On other hand, for all $t' > t$, one can use the dynamic programming principle \eqref{eq:DPP} and Condition \Href{ass:bound_G} to obtain that
	\begin{equation} \label{eq:DPP_alpha0}
		\vr(t,\xr) ~\ge~ \vr(t', \xr_{t \wedge \cdot}) - \int_t^{t'} G_s(\xr_{t \wedge \cdot}, 0) ds
		~\ge~
		\vr(t', \xr_{t \wedge \cdot}) - C(t'-t),
	\end{equation}
	and therefore
	$$
		\vr(t, \xr)  ~\ge~ \limsup_{t' \searrow t} \vr(t', \xr_{t \wedge \cdot}).
	$$
	It follows that $\vr \in \Cb(\Theta)$.
	Moreover, \eqref{eq:DPP_alpha0} implies   \eqref{eq:non_increas}.
	\endproof

	\vspace{2mm}
	
	For the next result, recall    the constant $C_0 > 0$   given in \Href{ass:bound_G}. Based on it, we define the (convex) function $\Gamma_0: \TD \to \R$ by 
	\begin{equation} \label{eq:def_Gamma}
		\Gamma_0(t, \xr) ~:=~ C_0 \xr_t^2.
	\end{equation}
 
	\begin{proposition}\label{prop : bar v - Gamma concave} 
		The functional $\vr-\Gamma_0$ is Dupire-concave, in the sense that,
		$$
			y \longmapsto (\vr - \Gamma_0) (t, \xr\otimes_t y) ~\mbox{is concave,}
		$$ 
		for all $t< T$ and $\xr \in \DT$.
	\end{proposition}
	\proof Let us take $t< T$, $\xr, \xr^1, \xr^2 \in \DT$ such that $\xr_s = \xr^1_s = \xr^2_s$ for $s \in [0,t)$,
	and $\xr_t = \theta \xr^1_t + (1-\theta) \xr^2_t$ for some $\theta \in [0,1]$.
	To prove the Dupire-concavity of $\vr - \Gamma_0$, it is enough to prove that
	\begin{equation} \label{eq:v_concave}
		\vr (t, \xr) - \Gamma_0(t, \xr)
		~\ge~
		\theta \big( \vr(t, \xr^1) -  \Gamma_0(t, \xr^1) \big)
		~+~
		(1-\theta) \big( \vr(t, \xr^2) - \Gamma_0(t, \xr^2) \big). 
	\end{equation}
	Without loss of generality, let us assume that 
	$$
		\xr^1_t < \xr_t < \xr^2_t.
	$$
	Let us now consider a standard Brownian motion $W^* = (W^*_s)_{s \ge 0}$ on some probability space $(\Om^*, \Fc^*, \P^*)$, and define the stopping time 
	$$
		\tau^* := \inf \big\{s \ge 0 ~: W^*_s = \xr^1_t - \xr_t ~\mbox{or}~ W^*_s = \xr^{2}_t - \xr_t \big\},
	$$
	as well as the martingale $X^*$ by
	$$
		X^*_s := W^*_{\tau^* \wedge \frac{s}{h-s}} \1_{\{s < h\}} +  W^*_{\tau^*} \1_{\{s \ge h\}},\;s\le T.
	$$
	Define $\P^*_h := \P^* \circ (X^*_{\cdot})^{-1}$. It is clear that $\P^*_h \in \Pc_0$ and that it satisfies
	\begin{equation} \label{eq:propty_Ph}
		\P^*_h \big[ \Xb^{t,\xr}_{t+h} = \xr^1_t \big] = \theta,
		~~
		\P^*_h \big[ \Xb^{t,\xr}_{t+h} = \xr^2_t \big] = 1- \theta
		~~\mbox{and}~~
		\P^*_h \big[ \Xb^{t,\xr}_s \in [\xr^1_t, \xr^2_t], ~s \in [t,t+h] \big] = 1,
	\end{equation}
	see e.g.~\cite[Exercise 8.13, Chapter 2.8]{karatzas1998brownian}.
	Then, it follows from the dynamic programming principle \eqref{eq:DPP}, together with the definition of $\Gamma_0$ in \eqref{eq:def_Gamma}, that
	\begin{eqnarray} \label{eq:v_concave_inter}
		\vr(t, \xr) - \Gamma_0(t, \xr) 
		&\ge&
		\E^{\P^*_h} \Big[ \vr(t+h, \Xb^{t,\xr}) - \Gamma_0(t+h, \Xb^{t,\xr}) - \int_t^{t+h} \big(G_s(\Xb^{t,\xr}, \alphab_s) - C_0 \alphab^2_s \big)ds \Big] \nonumber \\
		&\ge&
		\E^{\P^*_h} \Big[ \vr(t+h, \Xb^{t,\xr}) - \Gamma_0(t+h, \Xb^{t,\xr}) - Ch \Big] \nonumber \\
		&=&
		\theta \big( \vr(t+h, \xr^1_{t \wedge \cdot}) - \Gamma_0(t+h, \xr^1_{t \wedge \cdot}) \big) \nonumber\\
		&&
		+(1-\theta) \big(\vr(t+h, \xr^2_{t \wedge \cdot}) - \Gamma_0(t+h, \xr^2_{t \wedge \cdot}) \big)
		+C(h),
	\end{eqnarray}
	where 
	\begin{align*}
		C(h) 
		&:=
		 \E^{\P^*_h} \big[ \vr(t+h, \Xb^{t,\xr}) \big] - \theta  \vr(t+h, \xr^1_{t \wedge \cdot})  - (1 - \theta)  \vr(t+h, \xr^2_{t \wedge \cdot}) - Ch\\
		 &= \sum_{i=1}^{2}\E^{\P^*_h} \big[ (\vr(t+h, \Xb^{t,\xr}) -\vr(t+h, \xr^i_{t \wedge \cdot}) )\1_{\{\Xb^{t,\xr}_{t+h}=\xr^{i}_{t+h}\}}\big]  - Ch.
	\end{align*}
	We claim that $C(h) \to 0$ as $h \searrow 0$.
	Recalling that $\vr \in \Cb(\Theta)$, by (ii) of Proposition \ref{prop : continuity of J}, 
	then  \eqref{eq:v_concave} follows from  \eqref{eq:v_concave_inter} by letting $h \searrow 0$.

	\vspace{0.5em}
	
	To prove that $C(h) \to 0$, let us define
	$$
		D^{ \xr}_{t,h} 
		~:=~ 
		\big\{ \bar \xr \in \DT 
			~: \bar \xr_s = \xr_s ~\mbox{for}~s \in [0,t),
			~~\mbox{and}~~
			\bar \xr_s \in [\xr^1_t, \xr^2_t]~\mbox{for}~s \in [t, t+h]
		\big\},
	$$	
	and consider $\bar \xr^1, \bar \xr^2 \in D^{ \xr}_{t,h} $ such that $\bar \xr^1_{t+h} = \bar \xr^2_{t+h}$.
	Using Remark \ref{rem:optim_ctrl_bound} and recalling the definition of the sets $\widetilde \Pc_0^M$   and $D^M_{t+h}$   in \eqref{eq:def_DM},
	we deduce that one can find  $M>0$ such that
	\begin{eqnarray*}
		&&
		\big| \vr(t+h, \bar \xr^1) - \vr(t+h, \bar \xr^2) \big|\\
		&\le&
		\sup_{\P\in \widetilde \Pc_0^M} 
		\Big| \E^{\P} \Big[ 
			\Phi(\Xb^{t+h, \bar \xr^1}) - \Phi(\Xb^{t+h, \bar \xr^2})
			+
			\int_{t+h}^T G_s(\Xb^{t+h, \bar \xr^1}, \alphab^{t+h}_s) - G_s(\Xb^{t+h, \bar \xr^2}, \alphab^{t+h}_s) ds 
		\Big]
		\Big| \\
		&\le&
		\Big( \sup_{\bar \xr \in D^M_{t+h}} 
			\Big|
				\int_t^{t+h}  \!\!\! \lambda_{\Phi} (du, \bar \xr)
			\Big|	
		+ 
		\sup_{\P\in \widetilde \Pc_0^M} 
		\E^{\P}\Big[
		\sup_{\bar \xr \in D^M_{t+h}}
			\Big|
				\int_{t+h}^T \int_t^{t+h} \!\!\! \lambda_G(du, \bar \xr; s, \alphab^{t+h}_s) ds
			\Big|
		\Big]		\Big)
		~ | \xr^2_{t} -  \xr^1_{t}|.
	\end{eqnarray*}
	By Condition \Href{ass:concavity},
	  $\big| \vr(t+h, \bar \xr^1) - \vr(t+h, \bar \xr^2) \big| \to 0$ as $h \searrow 0$.
	Combining this with \eqref{eq:propty_Ph} implies that  $C(h) \to 0$ as $h \searrow 0$.
	\endproof

	\begin{remark}\label{rem: conc enveloppe pour cond term} 
		Fix $\P\in \Pc_{0}$ and $(t,\xr)\in [0,T)\x D([0,T])$. It follows from the dynamic programming principle in Proposition \ref{prop:equiv_dpp} that 
		\begin{equation*} 
			\vr(t,\xr)
			~\ge ~
			\E^{\P} \Big[
				\vr(T-h, \Xb^{t,\xr}_{(T-h) \wedge \cdot}) - \int_{t}^{T-h} G_{s}(\Xb^{t,\xr}, \alphab^t_{s}) ds
			\Big],
		\end{equation*}
		for all   $h\in (0,T-t)$. On the other hand,  \eqref{eq:non_increas} implies that $\liminf_{h\to 0} \vr(T-h, \Xb^{t,\xr}_{(T-h) \wedge \cdot})\ge \Phi(\Xb^{t,\xr})$. By Proposition \ref{prop : bar v - Gamma concave}, it follows that 
		\begin{equation*} 
			\vr(t,\xr)
			~\ge ~
			\E^{\P} \Big[
				\hat \Phi(T, \Xb^{t,\xr}) - \int_{t}^{T} G_{s}(\Xb^{t,\xr}, \alphab^t_{s}) ds
			\Big],
		\end{equation*}
		in which $\hat \Phi$ is the smallest function above $\Phi$ such that $y \longmapsto (\hat \Phi - \Gamma_0) (T, \xr\otimes_T y)$ is concave for all $\xr \in \DT$. It is given by  $(\Phi - \Gamma_0)^{\rm conc}+\Gamma_{0}$ in which $(\Phi - \Gamma_0)^{\rm conc}(\xr)$ is computed as the concave enveloppe of $y \longmapsto (\Phi - \Gamma_0) (T, \xr\otimes_T y)$ at $y=0$, see e.g.~\cite[Remark 3.8]{BLSZ18}. 
		This implies that 
		\begin{equation*} 
			\vr(t,\xr)
			~= ~
			\sup_{\P\in \Pc_{0}}\E^{\P} \Big[
				\hat \Phi(T, \Xb^{t,\xr}) - \int_{t}^{T} G_{s}(\Xb^{t,\xr}, \alphab^t_{s}) ds
			\Big],
		\end{equation*}
		and that any optimal rule $\P$ of problem \eqref{eq:dual_pb} satisfies $\P \big[\hat \Phi(T, \Xb^{t,\xr})\ne  \Phi(T, \Xb^{t,\xr}) \big]=0$.
	\end{remark}

\subsection{$\Cb^{0,1}(\Theta)$-regularity of $\vr$} 
 
	We are now in a position to prove that $\vr$ belongs to $\Cb^{0,1}(\Theta)$ and that the enveloppe principle \eqref{eq: gradient bar v} is valid. 
	The $\Cb^{0,1}(\Theta)$-regularity of $\vr$ allows us to apply the It\^{o} formula, 
	which combined with \eqref{eq: gradient bar v} and Proposition \ref{prop : continuity of J} will lead to the key property linking $\nabla_{\xr} \vr$, $A$ and $\hat \gamma$ in Proposition \ref{prop: characterization hat s}  below.

	\begin{proposition}\label{prop : regu et decroissance de bar v} 
		$\mathrm{(i)}$ Fix $t<T$, $\xr \in \DT$, and let $\Ph$ be an optimal solution to the control problem \eqref{eq:dual_pb} with initial condition $(t, \xr)$.
		Then, the map $(h,y)\in [0,T-t) \x \R \mapsto \vr(t+h, \xr\oplus_{t}y)$ is continuous, and $\vr$ is vertically differentiable at $(t,\xr)$ with vertical derivative given by 
		\begin{align}\label{eq: gradient bar v}
			\nabla_{\xr} \vr(t,\xr)= \nabla_{\xr} J(t,\xr;  \Ph). 
		\end{align}
		In particular, $\nabla_{\xr} \vr $ is locally bounded.
		Moreover, $\vr \in \Cb^{0,1}(\Theta)$.
	\end{proposition}
	\proof $\mathrm{(i)}$ First, the continuity of $(h, y) \mapsto \vr(t+h, \xr \oplus_t y)$ follows directly from the fact that $\vr \in \Cb(\Theta)$, see Proposition \ref{prop : continuity of J}.
	
	\vspace{0.5em}
	
	Let us now set
	$$
		\vp(y) ~:=~ \vr(t, \xr \oplus_t y) - J(t, \xr \oplus_t y; \Ph),~~\mbox{for all}~y \in \R.
	$$
	Then, $\vp$ achieves its minimum, equal to $0$, at $y=0$. 
	This implies that
	$$
		\limsup_{\eps \searrow 0} \frac{\vp(0) - \vp(-\eps)}{\eps}
		~\le~ 
		0
		~\le~
		\liminf_{\eps \searrow 0} \frac{ \vp(\eps) - \vp(0)}{\eps},	
	$$
	and hence
	$$
		\limsup_{\eps \searrow 0} \frac{\vr(t, \xr)- \vr(t, \xr \oplus_t(-\eps))}{\eps}
		~\le~ 
		\nabla_{\xr} J(t,\xr; \Ph)
		~\le~
		\liminf_{\eps \searrow 0} \frac{ \vr(t, \xr \oplus_t \eps) - \vr(t, \xr)}{\eps}.
	$$
	On the other hand, $y\mapsto \vr(t,\xr\oplus_{t}y)-\Gamma_0(t, \xr\oplus_{t}y)$ is concave, and $y \mapsto \Gamma_0(t, \xr \oplus_t y)$ is differentiable, which implies that
	$$
		\liminf_{\eps \searrow 0} \frac{\vr(t, \xr)- \vr(t, \xr \oplus_t(-\eps))}{\eps}
		~\ge~
		\limsup_{\eps \searrow 0} \frac{ \vr(t, \xr \oplus_t \eps) - \vr(t, \xr)}{\eps}.
	$$
	This proves \eqref{eq: gradient bar v}.
	The fact that   $\nabla_{\xr} \vr$ is locally bounded is   an immediate consequence of  \eqref{eq:v_Lip}  and Remark \ref{rem:optim_ctrl_bound}.

	\vspace{0.5em}

	To prove that $\nabla_x \vr \in \Cb(\Theta)$, 
	let us consider $(t^{n},\xr^{n})_{n\ge 1}\subset \Theta$   such that $t^{n}\ge t$, $(t^{n},\xr^{n})_{n \ge 1}$ converges to $(t,\xr)$ on $(\Theta, d)$, and $y^{n}\to y \in \R$.
	 Notice that $(\vr(t_{n},\xr^{n}+\1_{\{t_{n}\}}\cdot)-\Gamma_{0} (t_{n}, \xr^{n}+\1_{\{t_{n}\}}\cdot))_{n\ge 1}$ is a sequence of concave functions which converges to the concave function $\vr(t,\xr\oplus_{t} \cdot)-\Gamma_0(t, \xr\oplus_{t} \cdot)$.
	 Then any limit of the corresponding sequence of gradients $(\partial_{y} \vr(t_{n},\xr^{n}+\1_{\{t_{n}\}}\cdot)-\partial_{y}\Gamma_0 (t_{n}, \xr^{n}+\1_{\{t_{n}\}}\cdot))_{n\ge 1}$, computed at $y=0$, converges to an element of the super-differential of $y\mapsto  \vr(t,\xr\oplus_{t}y)-\Gamma_{0}({t}, \xr\oplus_{t}y)$, computed at $y=0$.
	 As $y\mapsto  \vr(t,\xr\oplus_{t}y)-\Gamma_0({t}, \xr\oplus_{t}y)$ is differentiable, its super-differential at $0$ is unique and equal to $\nabla_{\xr}  \vr(t,\xr)- \partial_{y}\Gamma_0({t}, \xr\oplus_{t}y)_{|y=0}$.
	\endproof

	\vspace{0.5em}
 
	\begin{proposition}\label{prop : ito surbar v}  
		Fix $x_0 \in \R$ and let $\Ph$ be an optimal solution to the control problem \eqref{eq:dual_pb} with initial condition $(0, x_0)$.
		Then
		$$
			\vr(t, X)
			~=~
			\vr(0, x_0)
			+
			\int_0^t \nabla_{\xr}  \vr(s, X)d X_s
			+
			\int_0^t G_s( X,\alpha_s)ds,
			~~t \in [0,T],
			~~
			\Ph \mbox{-a.s.}
		$$
	\end{proposition}
 	\proof 
	By Proposition \ref{prop:equiv_dpp} and Lemma \ref{lemm:optimal_ctrl_cond}, the process
	$$
		 \vr(t, X) - \int_0^t  G_s( X,\alpha_s)ds, ~t \in [0,T], ~~\mbox{is a}~\Ph \mbox{-martingale}.
	$$
	Then, it is enough to appeal to the It\^o's formula in Proposition \ref{prop: Dupire Tanaka corrige} (together with \eqref{eq:non_increas}, Proposition \ref{prop : bar v - Gamma concave} and Proposition \ref{prop : regu et decroissance de bar v}).
	\endproof


\section{Proof of Theorem \ref{thm:dualty}}

	All over this section, we fix $x_0 \in \R$ as the initial condition of the dual control problem, and we assume that there exists an optimal control rule  to the control problem \eqref{eq:dual_pb} with initial condition $(0,x_0)$. The assumptions of Theorem \ref{thm:dualty} are in force. 

\subsection{The key ingredient}
  
	We start by providing the key ingredient for the proof of Theorem \ref{thm:dualty}, which relies $\nabla_{\xr} \vr$, $A$ and $\hat \gamma$ in Proposition \ref{prop: characterization hat s} below (recall the notations in \eqref{eq:def_At}-\eqref{eq:def_gamma}). For this, we first need to check that a dual optimizer can be chosen to be extremal for $X$.

	\begin{lemma}\label{lem : MRP} There exists an optimal control rule $\Ph \in \Pc^*_{0, x_0}$ which is also an extremal point
		of the space of all martingale measures, i.e.~probability measures on $\Om$ under which $X$ is a martingale.
		In particular, any $(\Ph, \F)$-martingale ${\cal M} = ({\cal M}_t)_{t \in [0,T]}$ on $\Om$ has the representation
		\begin{equation} \label{eq:mgt_repre}
			{\cal M}_t = {\cal M}_0 + \int_0^t \beta_s dX_s, ~t \in [0,T], ~~\Ph\mbox{-a.s.}
		\end{equation}
		for some $\F$-predictable process $\beta$ satisfying $\int_0^T \beta^2_s d\langle X \rangle_s < \infty$, $\Ph$-a.s.
	\end{lemma}
	\proof Let us denote by $\Ec_{0, x_0}$ the collection of all extreme points of the space $\Pc_{0,x_0}$,
	and let $\Pt \in \Pc^*_{0,x_0}$ be an optimal control rule for the control problem \eqref{eq:dual_pb} with initial condition $(0, x_0)$.
	Then, by  \cite[Corollary (3.3) and Remark p. 111]{jacod1977etude}, 
	there exists a measure $\mu$ on $\Pc_{0}$ such that $\Pt[A]= \int_{\Ec_{0,x_{0}}} \P[A]\mu(d\P)$ for all Borel subset $A$ of $\Om$.
	As the reward function $J(0, x_0; \P)$ is linear in $\P$, 
	it follows that, for $\mu$-a.e. $\P$ in $\Ec_{0, x_0}$, $\P$ is an optimal solution in $\Pc^*_{0, x_0}$.
	Therefore, there exists an optimal control rule $\Ph \in \Pc^*_{0, x_0}$, which also belongs to $\Ec_{0, x_0}$.
	\endproof

	\vspace{0.5em}

	Let us now fix  an optimal control rule $\Ph$ as in Lemma \ref{lem : MRP}, i.e.~which is also an extreme martingale measure.
	 Recall that $\gammah$ is defined in \eqref{eq:def_gamma},
	$A$  in \eqref{eq:def_At}, and that  $B^{\P}$ is the dual predictable projection of $A$ w.r.t. $(\P, \F)$ (see Remark \ref{rem:dual_proj}).

	\begin{proposition}\label{prop: characterization hat s} 
		The random variable $A_T$ is $\Ph$-integrable.
		Moreover,   $\nabla_{\xr} \vr(0, x_0) = \E^{\Ph}[ A_T]$ $ {=\E^{\Ph}[B^{\Ph}_T]}$,
		and the following martingale representation holds:
		\begin{equation} \label{eq:prop_AT}
			 {B^{\Ph}_T} ~=~ \nabla_{\xr} \vr(0, x_0) + \int_0^T  \gammah_s dX_s, ~\Ph\mbox{-a.s.},
		\end{equation}
		where $\big(\int_0^t \gammah_s dX_s \big)_{t \in [0,T]}$ is a $\Ph$-martingale.
	\end{proposition}
	\proof 
	Let us write $\Eh$ for $\E^{\Ph}$ for ease of notations.
	The integrability of $A_T$ under $\Ph$ follows directly from \Href{ass:Phi_Frechet_diff} and \Href{ass:G_Frechet_diff}.
	Moreover, \eqref{eq: AT -At}, \eqref{eq:DJ}, \eqref{eq: gradient bar v}  and direct computations imply that
	$$
		\nabla_{\xr} \vr(0, x_0) 
		~=~
		\nabla_{\xr} J(0, x_0; \Ph) 
		~=~
		\Eh \Big[ \int_0^{T}  \lambda_{\Phi}(ds,  X)- \int_0^{T}  \int_{[0,s]} \lambda_{G}(du, X; s, \alpha_s)ds \Big] 
		~=~
		\Eh[A_T].
	$$
	
	 {Since $\widehat \E[B^{\Ph}_T]=\widehat \E[A_{T}]$},  the martingale representation \eqref{eq:mgt_repre} implies 
		$$
		{B^{\Ph}_T}
		~=~
		\nabla_{\xr} \vr(0, x_0) 
		+ 
		\int_0^T \beta_s dX_s,
		~~\Ph\mbox{-a.s.},
	$$
	for some $\F$-predictable process $\beta$ such that $(\int_0^t \beta_s dX_s)_{t \in [0,T]}$ is a $\Ph$-martingale, and it just remains  to prove that $\beta = \gammah $, $d\P \x d \langle X \rangle$-a.e. 
	
	\vspace{0.5em}
	
	Let us denote by $\Ac_{\infty}$ the collection of all bounded $\R$-valued $\F$-predictable processes.
	Given $\delta \in \Ac_{\infty}$, define $Z^{\delta} := \int_{0}^{\cdot}\delta_s dX_s$.
	Then, for $\eps>0$  and $\theta \in [- \eps, \eps]$, $\Ph \circ (X + \theta Z^{\delta})^{-1}$ is a weak control rule,
	and the map
	$$
		\theta \in \R 
		~\longmapsto~
		\Eh \Big[ 
			\Phi( X+ \theta Z^{\delta})  
			-
			\int_{0}^{T}G_{s} \big( X+ \theta Z^{\delta},  \alpha_{s} (1 +\theta \delta_{s}) \big)ds
		\Big]
	$$
	achieves its minimum at $\theta = 0$.
	By  standard arguments in calculus of variations, together with \Href{ass:Phi_Frechet_diff}, \Href{ass:G_Frechet_diff} and the arguments in Remark \ref{rem: lin growth when bounded frechet derivative}, one deduces that 
	\begin{equation}\label{eq:mgt_repre_inter}
		0
		~=~
		\Eh \Big[ 
			\int_{0}^{T} Z^{\delta}_{s} \lambda_{\Phi}(ds, X)
			~-~
			\int_{0}^{T} \Big( \int_{0}^{s}  Z^{\delta}_{u} \lambda_{G}(du, X; s, \alpha_{s})+\delta_{s} \alpha_s \partial_{a}G(s, X, \alpha_{s}) \Big)ds 
		\Big].
	\end{equation}
	Let us  {recall that}
	$$
		\widehat \lambda_{G}(I) ~=~ \int_0^T  {\int_0^T}  \1_{[0,s] \cap I}(u) \lambda_{G}(du, X; s, \alpha_{s})ds,
		~\mbox{for all Borel subsets}~I \subset [0,T],
	$$
	so that 
	$$
		\int_{0}^{T}    \int_{0}^{s}  Z^{\delta}_{u} \lambda_{G}(du, X; s, \alpha_{s})ds
		~=~
		\int_{0}^{T}   \int_{0}^{T}  Z^{\delta}_{u}\1_{[0,s]}(u) \lambda_{G}(du, X; s, \alpha_{s})ds
		~=~
		\int_{0}^{T}  Z^{\delta}_{u}\widehat \lambda_{G}(du),
	$$
	and let us   denote by $  \lambda^{\circ}_{\Phi}{(\cdot,X)}$ (resp. $\widehat \lambda^{\circ}_G$) the dual predictable projection of $\lambda_{\Phi}{(\cdot,X)}$ (resp. $\widehat \lambda_G$) w.r.t. $(\F, \Ph)$,
	{so that $B^{\Ph}_T =  \int_0^T  ( \lambda^{\circ}_{\Phi}(ds, X) - \widehat \lambda^{\circ}_G (ds) )$.}
	Then, it follows from  \eqref{eq:mgt_repre_inter} and \eqref{eq:def_At} that
	\begin{eqnarray*}
		&&
		\Eh \Big[ \int_0^T \delta_s \alpha_s \partial_a G(s, X, \alpha_s) ds \Big]
		~=~
		\Eh \Big[ \int_0^T Z^{\delta}_s \Big( \lambda_{\Phi}(ds, X) - \widehat \lambda_G (ds) \Big) \Big] \\
		&=&
		\Eh \Big[ \int_0^T Z^{\delta}_s \Big( \lambda^{\circ}_{\Phi}(ds, X) - \widehat \lambda^{\circ}_G (ds) \Big) \Big]	
		~=~	
		\Eh \Big[ Z^{\delta}_T \int_0^T  \Big( \lambda^{\circ}_{\Phi}(ds, X) - \widehat \lambda^{\circ}_G (ds) \Big) \Big]\\
		&=&
		\Eh \big[ Z^{\delta}_T {B^{\Ph}_T} \big]
		~=~
		\Eh \Big[ \int_0^T \delta_s \beta_s \alpha^2_s ds \Big].
	\end{eqnarray*}
	 By the arbitrariness of $\delta$, the definition of $\gammah$ in \eqref{eq:def_gamma}, and since  $\alpha$ is the density of $\langle X \rangle$, this proves that 
	$$
		 \beta= \frac{\partial_a G(\cdot, X, \alpha)}{\alpha} =  \gammah,~d\Ph \x d \langle X \rangle \mbox{-a.e.}
	$$
	\endproof

\subsection{Construction of the hedging strategy}\label{sec: construction of the hedging strat}

	We now have all the required ingredients for the construction of a perfect hedging strategy based on $\hat \gamma$ and the dual predictable projection $B^{\Ph}$ of $A$, under an (extreme) optimal measure $\Ph$.

	\vs2
	
	{\bf Proof of Theorem \ref{thm:dualty}.} Let us consider the extreme martingale measure and   optimal control rule $\Ph$ as in Lemma \ref{lem : MRP}.  Define
	$$
		v_0 := \vr(0,x_0),
		~~
		V_t ~:=~ \vr(t, X),
		~~
		y_0 := \nabla_{\xr} \vr(0,x_0),
		~~
		Y_t := \nabla_{\xr} \vr(t, X),
		~~\mbox{for all}~t \in [0,T].
	$$
	Then, it follows from Proposition \ref{prop : ito surbar v}, \eqref{eq:structure} and \eqref{eq:def_gamma} that
	$$
		V_t = v_0 + \int_0^t Y_s dX_s + \int_0^t F_s(X, \gammah_s) ds, 
		~~t \in [0,T],
		~~\Ph \mbox{-a.s.}
	$$
	Recall that $A = (A_t)_{t \in [0,T]}$ is the bounded variation process defined in  \eqref{eq:def_At}  and  that
	$B^{\Ph} = (B^{\Ph}_t)_{t \in [0,T]}$ is the dual predictable projection of $A$ w.r.t. $(\F, \Ph)$.
	We claim that
	\begin{equation} \label{eq:proof_thm_claim}
		Y_t  ~=~ \widehat Y_t ~:=~ y_0 + \int_0^t \gammah_s dX_s - B^{\Ph}_{{t-}}, ~~t \in [0,T], ~~\Ph\mbox{-a.s.}
	\end{equation}
	Since $\vr(T,X)=\Phi(X)$, it is then clear that the probability measure $\Ph$, together with $(v_0, y_0)$ and {$(V_t, Y_t, \gammah_t, B^{\Ph}_{t-})_{t \in [0,T]}$} defined above, provides a solution to the initial hedging problem in Definition \ref{def:hedge_pb}.	
	\vspace{0.5em}
	
	Thus, it remains to prove the claim \eqref{eq:proof_thm_claim}.
	Recall that $\nabla_{\xr} \vr \in \Cb(\Theta)$ and is non-anticipative, which
	 implies that $Y:= \nabla_{\xr} \vr(\cdot, X)$ is $\F$-adapted and right-continuous, and hence is $\F$-progressively measurable, or equivalently is $\F$-predictable.
	At the same time $\widehat Y$ is also $\F$-predictable by its definition in \eqref{eq:proof_thm_claim}.
	
	\vspace{0.5em}
	
	We next consider a $\F$-stopping time $\tau$ taking value in $[0,T]$, and
	let $(\Ph_{\om})_{\om \in \Om}$ be a r.c.p.d. of $\Ph$ knowing $\Fc_{\tau}$.
	By Lemma \ref{lemm:optimal_ctrl_cond}, for $\Ph$-a.e.~$\om$, $\Ph_{\om}$ is an optimal control rule for the control problem \eqref{eq:dual_pb_equiv} with initial condition $(\tau(\om), X_{\tau(\om)\wedge \cdot}(\om))$.
	As in Proposition \ref{prop : continuity of J}, 
	it follows by direct computations that, for $\Ph$-a.e. $\om \in \Om$,
	\begin{eqnarray*} 
		&&
		\nabla_{\xr} \vr(\tau(\om), X(\om))
		~=~
		\nabla_{\xr}J (\tau(\om), X(\om); \Ph_{\om})  \\
		&=&
		\E^{\Ph_{\om}} \Big[ \int_{\tau(\om)}^T \lambda_{\Phi}(du, X) - \int_{\tau(\om)}^T \int_{\tau(\om)}^s \lambda_G(du, X; s, \alpha_s) ds \Big] 
		~=~
		\E^{\Ph_{\om}} \big[ A_T - A_{\tau(\om){-}} \big],
	\end{eqnarray*}
	{in which the last equality follows from \eqref{eq:def_At}.}
	Then by \cite[Theorem VI.75]{dellacherie1982probabilities}, it follows that
	\begin{equation} \label{eq:prop_A}
		\nabla_{\xr} \vr(\tau, X)
		~=~
		\E^{\Ph} \Big[ A_T - A_{\tau{-}} \Big| \Fc_{\tau} \Big]
		~=~
		\E^{\Ph} \Big[ B^{\Ph}_T - B^{\Ph}_{\tau{-}} \Big| \Fc_{\tau} \Big]
		,
		~~\Ph\mbox{-a.s.~~}
	\end{equation}

	Recall that $A_T$   is $\Ph$-integrable, see Proposition \ref{prop: characterization hat s}.
	Using \eqref{eq:prop_A} and then Proposition \ref{prop: characterization hat s}, it follows that, for all $\F$-stopping times $\tau$,
	\begin{eqnarray*}
		Y_{\tau} 
		=
		\nabla_{\xr} \vr(\tau, X)
		\!\!\!\!&=&\!\!\!\!
		\E^{\Ph} \big[ B^{\Ph}_T -  {B^{\Ph}_{\tau{-}}}  \big| \Fc_{\tau} \big] 
		=
		\E^{\Ph}\big[  B^{\Ph}_T \big| \Fc_{\tau} \big]  -B^{\Ph}_{\tau{-}}
		=
		\E^{\Ph} \Big[ y_0 + \int_0^T \gammah_s dX_s \Big| \Fc_{\tau} \Big]  -B^{\Ph}_{\tau{-}} \\
		\!\!\!&=&\!\!\!
		y_0 + \int_0^{\tau} \gammah_s dX_s  -B^{\Ph}_{\tau{-}}
		~=~
		\widehat Y_{\tau}, ~~\Ph \mbox{-a.s.}
	\end{eqnarray*}
	Since  $Y$ and $\widehat Y$ are both $\F$-predictable processes, 
	we can use the predictable section theorem (see e.g. \cite[Theorem IV.86]{dellacherieMeyer}) to conclude that $Y_{\cdot} = \widehat Y_{\cdot}$, $\Ph$-a.s., which is the claim in \eqref{eq:proof_thm_claim}.
	\endproof


\section{Sufficient conditions for the existence of a  solution to the dual optimal control problem}
\label{sec: existence}

	The stochastic optimal control problem \eqref{eq:dual_pb} does not match the usual assumptions that ensure existence of a weak optimal control rule, because the penalty term is only of quadratic growth, compare (H3) with \cite[Conditions (2.1) and (3.5)]{haussmann1990existence}.
	In this section, we  provide two different sufficient conditions. 
	The first one is a concavity condition on the reward function, which is quite classical and actually ensures existence of a strong optimal control rule.
	The second one is inspired from the PDE estimates in \cite{BLSZ18} and seems to be new in the optimal control literature.

\subsection{Existence of a strong optimal control rule}

	We start with a first result on the existence of a strong solution whenever sufficient concavity on the coefficients is assumed.

	\begin{proposition}\label{prop : existence strong} 
		Assume that $\Phi$ and $-G_{t}$ are concave, for all $t\le T$.
		Then, for all $(t,\xr) \in \TD$, there exists an optimal control rule $\Ph \in \Pcb^*_{t,\xr}$ for the dual problem \eqref{eq:dual_pb}, 
		which is also a strong control rule, i.e. $\Ph \in \Pcb^S_{t,\xr}$.
	\end{proposition}
	
	Let us comment on situations in which the above Proposition can be applied, before providing its proof.

	\begin{remark}
		In the context of Example \ref{exam:BoLoZo},  with the computation in Remark \ref{rem:G},
		it is clear that $(\xr,a) \mapsto -G_t(\xr, a)$ is concave whenever $f(\cdot) \equiv f_0$ for some constant $f_0 > 0$ and $x \mapsto \sigma_0(t,x)$ is an affine function.
	\end{remark}
	
	\begin{remark}\label{rem: strong existence}
		In the context of Example \ref{exam:BoLoZo}, assume that $ f(\cdot) \equiv f_0$ and $\sigma_{0}(\cdot)\equiv \sigma_{0}$ for some constants $f_0,\sigma_{0} > 0$, then 
		$$
		\E^{\P} \Big[\Phi\big( \Xb^{t,\xr} \big) - \int_{t}^{T}G_{s} \big( \Xb^{t,\xr}, \alphab^t_s \big) ds \Big]
		=\E^{\P} \Big[\Phi_{f_{0}}\big( \Xb^{t,\xr} \big)- \frac1{2f_{0}}\int_{t}^{T} \big(-2\alphab^t_s \sigma_{0}+\sigma_{0}^{2} \big) ds \Big]+\frac1{2f_{0}}\xr_{t}^{2},
		$$
		in which $\Phi_{f_{0}}\big( \xr' \big):=\Phi\big( \xr' \big)-\frac1{2f_{0}} ( \xr'_{T} )^{2}$ for $\xr'\in D([0,T])$. 
		By the arguments in the proof of  Proposition \ref{prop : existence strong} below, it suffices to assume that $\Phi_{f_{0}}$ is concave in order to prove the same existence results as in Proposition \ref{prop : existence strong}.
		On the other hand, it follows from Remark \ref{rem: conc enveloppe pour cond term} that 
		\begin{align*}
		\vr(t,\xr)&=\sup_{\P \in \Pc_0}  \E^{\P} \Big[\Phi^{\rm conc}_{f_{0}}\big( \Xb^{t,\xr}\big)+\frac1{2f_{0}} |\Xb^{t,\xr}_{T}|^{2}- \frac1{2f_{0}}\int_{t}^{T} \big(\alphab^t_s -\sigma_{0} \big)^{2} ds \Big]
		\\
		&=\sup_{\P \in \Pc_0}  \E^{\P} \Big[\Phi^{\rm conc}_{f_{0}}\big( \Xb^{t,\xr}\big)- \frac1{2f_{0}}\int_{t}^{T} \big(-2\alphab^t_s \sigma_{0}+\sigma_{0}^{2} \big) ds \Big]+\frac1{2f_{0}}\xr_{t}^{2}
		\end{align*}
		in which $\Phi^{\rm conc}_{f_{0}}$ is the concave enveloppe of $x\mapsto \Phi_{f_{0}}(x\1_{\{T\}})$ and $\hat \Phi:\xr'\mapsto \Phi^{\rm conc}_{f_{0}}(\xr')+\frac1{2f_{0}}|\xr'_{T}|^{2}$. Hence, the second equation above  implies that (strong) existence holds with the terminal payoff $\hat \Phi$, and that it leads to the same value function. This is consistent with \cite{BLSZ18} in which the payoff that is perfectly hedged is 
	$\hat \Phi$.
	\end{remark}
	
	\begin{remark} Note that the PDE arguments used in \cite{BLSZ18} also allow to provide the existence of a strong solution in a Markovian setting without the above strong concavity conditions, but under rather restrictive regularity assumptions. We refer to Remark \ref{rem: eps conca} below for a presentation of the results of \cite{BLSZ18} in a more specific model (chosen so as to alleviate the notations and to match with our benchmark Example \ref{exam:BoLoZo}).
	\end{remark}
	
{\bf Proof of Proposition \ref{prop : existence strong}.}
	Without loss of generality, we can restrict to the case with initial condition $(0, x_0)$.
	Let us first introduce an equivalent strong formulation to the control problem \eqref{eq:dual_pb} with initial condition $(0, x_0)$.
	Recall that $\Om = C([0,T])$ is the canonical space of continuous paths {with canonical process $X$. Let us denote by $\P_0$ the Wiener measure, under which $X$} is a standard Brownian motion.
	On the filtered probability space $(\Om, \Fc, \F, \P_0)$, let us denote
	$$
		\Ac_2 
		~:=~ 
		\Big\{ 
			\delta = (\delta_t)_{t \in [0,T]} 
			~: 
			\delta ~\mbox{is}~\F \mbox{-predictable and}~ \E^{\P_0} \Big[ \int_0^T \delta^2_t dt \Big] < \infty
		\Big\},
	$$
	and then define
	$$
		X^{\delta}_t ~:=~ x_0 + \int_0^t \delta_s {dX_s},~\P_0\mbox{-a.s. for every}~ \delta \in \Ac_2.
	$$
	Recall the set $\Pc_0^S$ of  strong control rules defined in Definition \ref{def : strook varadhan} and the set $\Pc^S_{0, x_0}$ defined in Remark \ref{rem:Pc}. One has
	$$
		\Pc^S_{0, x_0} 
		~=~
		\big\{ \P_0 \circ \big( X^{\delta}_{\cdot} \big)^{-1} ~: \delta \in \Ac_2 \big\}.
	$$
	Moreover, one has the equivalence result between the weak formulation (in \eqref{eq:dual_pb} or \eqref{eq:dual_pb_equiv}) and the strong formulation of the control problem (see e.g. \cite[Theorem 4.5]{el2013capacities}):
	\begin{equation} \label{eq:dual_pb_strong}
		\vr(0, x_0) 
		~=~
		\sup_{\P \in \Pc^S_{0, x_0}} \overline J(0, x_0; \P)
		~=~
		\sup_{\delta \in \Ac_2}
		\E^{\P} \Big[
			\Phi(X^{\delta}) - \int_0^T G_s(X^{\delta}, \delta_s) ds
		\Big].
	\end{equation}

	Let us now complete the proof  of Proposition \ref{prop : existence strong}, which is in fact an immediate consequence of Komlos' lemma. 

	\vspace{0.5em}

	Recall that, by Remark \ref{rem:optim_ctrl_bound}, there is a constant $C>0$ such that the optimum is achieved in the class of controls $\delta$ satisfying the integrability condition $\E^{\P_0}[\int_0^{T}\delta^{2}_{s}ds]\le C$. 
	Let $(\delta^{n})_{n\ge 1}$ be a maximizing sequence satisfying the latter integrability condition. 
	It follows from Komlos' lemma that, up to passing to convex combinations, one can extract a subsequence, $(\tilde \delta^{n})_{n\ge 1}$, that converges $dt\times d\P_{0}$-a.e.~to some $\F$-predictable process $\delta$ satisfying 
	$$
		\E^{\P_0} \Big[\int_{0}^{T} \delta^{2}_{s}ds \Big]\le C,
	$$
	by convexity of $y\mapsto y^{2}$ and Fatou's Lemma. 
	Next, noticing that 
	$
		 \E^{\P_0}[\int_0^{T}|\delta_{s}-\tilde  \delta^{n}_{s}|^{2}ds]\le 4 C, 
	$
	it follows by Fatou's Lemma again that
	$$
		\E^{\P_{0}} \Big[\int_0^{T}|\delta_{s}-\tilde  \delta^{n}_{s}|^{2}ds \Big] 
		~=~
		\E^{\P_0}
		\big[
			\Ninfty{ X^{\delta}-  X^{ \tilde \delta^{n}}}^{2}
		\big]
		~\longrightarrow~ 0,
		~~~\mbox{as}~~
		n \longrightarrow \infty.
	$$
	By concavity of $\Phi$ and $-G$, $(\tilde \delta^{n})_{n\ge 1}$ is still a maximizing sequence
	for the optimization problem in the r.h.s. of \eqref{eq:dual_pb_strong}.
	By \Href{ass:Phi_Frechet_diff}, Remark \ref{rem: lin growth when bounded frechet derivative} and \Href{ass:bound_G}, there exists some constant $C>0$ such that $|\Phi(\xr)| \le C(1 + \|\xr\|)$ and $|G_t(\xr, a) | \le C(1 + a^2)$ for all $(t,\xr, a)\in [0,T]\x D([0,T])\x \R$.
	Then, it is enough to apply Fatou's Lemma again to prove that $\delta$ is an optimal control process to the control problem in the r.h.s. of \eqref{eq:dual_pb_strong}.
\endproof

\subsection{Existence of a weak optimal control rule}
	
	Let us now consider a path-dependent extension of the model considered  in Example \ref{exam:BoLoZo} and Remark \ref{rem:G},
	where\footnote{More general situations could be tackled with the same tools, we choose this formulation for sake of simplicity of the exposition.}
	\begin{equation} \label{eq:G_example}
		G(t,\xr, a) 
		~=~
		\frac12\gamma_2(t,\xr) a^2 - \gamma_1(t,\xr) a + \gamma_0(t,\xr),
	\end{equation}
	for some positive bounded continuous functionals $\gamma_i : \Theta \to \R_+$, $i=0,1,2$.
	Assume in addition that, for some modulus of continuity  $\rho: \R_+ \to \R_+$ and some  constant $C>0$,
	\begin{equation} \label{eq:Regularity_gamma}
		| \gamma_1(t, \xr) - \gamma_1(t', \xr') | ~\le~ \rho(|t-t'|) + C \|\xr_{t \wedge \cdot} - \xr'_{t' \wedge \cdot}\|.
	\end{equation}
	In this context, one can apply exactly the same arguments in Proposition \ref{prop : bar v - Gamma concave} to deduce that $\vr - \Gammab_0$ is Dupire-concave, 
	i.e., for all $t \ge 0$, $y \mapsto (\vr - \Gammab_0)(t, \xr \oplus_t y)$ is concave, with 
	$$
		\Gammab_0(t, \xr) 
		~:=
		\int_0^{\xr_t} \int_0^{y_1} \gamma_2 \big(t, \xr \oplus_t (y_2 - \xr_t) \big) dy_2 dy_1.
	$$
	We will nevertheless assume a slightly more restrictive condition, namely that		
	\begin{align}\label{eq: hyp eps conc}
		y\in \R \mapsto (\vr - \Gammab_{\eps_0})(t, \xr \oplus_t y)\;\;\mbox{ is concave for all $(t,\xr) \in [0,T]\x D([0,T])$.}
	\end{align}
	with
	$$
		\Gammab_{\eps_0}(t,\xr) ~:=~ \Gammab_0(t,\xr) - \eps_0  \xr_t^2,
	$$
	for some $\eps_0>0$. This condition will be further discussed in Remark \ref{rem: eps conca} below.

	\begin{theorem}
		In the context of \eqref{eq:G_example}, assume that \eqref{eq: hyp eps conc} holds for some $\eps_0>0$.
		Then for every initial condition $(t,\xr)$, there exists an optimal control rule $\Ph$ for the problem \eqref{eq:dual_pb}.
	\end{theorem}
	\proof Without loss of generality, we can restrict to the case with initial condition $(0, x_0)$.
	Let $(\P^n)_{n \ge 1}$ be a maximizing sequence for the control problem \eqref{eq:dual_pb_equiv}.
	By the boundedness of $\gamma_2$, together with Conditions \Href{ass:Phi_Frechet_diff}, \Href{ass:bound_G} and Remark \ref{rem: lin growth when bounded frechet derivative}, one has
	\begin{equation} \label{eq:bounded_Phi_X2}
		|\Phi(\xr)|  ~\le~ C(1 + \|\xr\|),
		~~~
		|\Gammab_0(t, \xr)| \le C(1+\|\xr_{t \wedge \cdot}\|^2)
		~~\mbox{and}~~
		\sup_{n \ge 1} \E^{\P_n} \big[ \sup_{0 \le t \le T} X_t^2 \big] \le C,
	\end{equation}
	 for some constant $C> 0$.
	Fix $M > 0$, $\tau_M := \inf \{t \ge 0~: |X_t| \ge M\}$ and let $\Tc_M$ denote the collection of all $\F$-stopping times dominated by $\tau_M$. We claim that, for all $M>0$,
	\begin{equation} \label{eq:claim_tight}
		\lim_{\theta \searrow 0} 
		\delta_M(\theta) =0,
		~~~\mbox{with}~~
		\delta_M(\theta)
		:= \limsup_{n \to \infty} 
		\sup_{\sigma, \tau \in \Tc_M, \sigma \le \tau \le \sigma + \theta} 
		\E^{\P_n} \big[ \big| X_{\tau} - X_{\sigma} \big|^2 \big].
	\end{equation}
	Then, using Aldous' Criterion (see e.g. \cite[Theorem VI.4.5]{jacod2013limit}) together with \cite[Proposition VI.3.26]{jacod2013limit} and \eqref{eq:bounded_Phi_X2}, one obtains that
	\begin{equation} \label{eq:tightness}
		\mbox{the sequence}~
		\Big( \P^n \circ \Big(X_{\cdot}, \int_0^{\cdot} \alpha^2_s ds \Big)^{-1} \Big)_{n \ge 1}
		~\mbox{is tight}.
	\end{equation}
	Therefore, up to possibly passing along a subsequence, there exists some $\Ph \in \Pc_{0,x_0}$ such that 
	\begin{equation} \label{eq:cvg_Pn}
		\P^n \circ \Big(X_{\cdot}, \int_0^{\cdot} \alpha^2_s ds \Big)^{-1}
		~\longrightarrow~
		\Ph \circ \Big(X_{\cdot}, \int_0^{\cdot} \alpha^2_s ds \Big)^{-1}.
	\end{equation}
	It follows that
	\begin{equation} \label{eq:cvg1}
		\lim_{n \to \infty} \E^{\P^n} \big[ \Phi(X) \big] = \E^{\Ph} \big[ \Phi(X) \big],
		~~~~~
		\lim_{n \to \infty} \E^{\P^n} \Big[ \int_0^T - \gamma_0(t, X)dt \Big]
		=
		\E^{\Ph} \Big[ \int_0^T - \gamma_0(t, X)dt \Big],
	\end{equation}
	and, since $\Gammab_0$ is bounded from below, that 
	\begin{eqnarray} \label{eq:cvg2}
		&&
		\limsup_{n \to \infty} \E^{\P^n} \Big[ \int_0^T - \frac12 \gamma_2(t, X) \alpha^2_t dt \Big]
		~=~
		\limsup_{n \to \infty} \E^{\P^n} \big[ \Gammab_0(0, x_0) - \Gammab_0(T,X) \Big] \nonumber \\
		&\le&
		\E^{\Ph} \big[ \Gammab_0(0, x_0) - \Gammab_0(T,X) \Big]
		~=~
		\E^{\Ph} \Big[ \int_0^T - \frac12 \gamma_2(t, X) \alpha^2_t dt \Big].
	\end{eqnarray}
	Moreover, recalling \eqref{eq:Regularity_gamma} and setting 
	\begin{eqnarray*}
		\delta(n, \eps) 
		&:=& 
		C \E^{\P^n} \Big[ \int_0^{\eps} \!\! \alpha_t dt \Big] 
		+ 
		\sqrt{ \E^{\P^n} \Big[ \int_{\eps}^T \alpha_t^2 dt \Big]} \sqrt{  \rho(\eps) + \E^{\P^n} \Big[\!\! \int_{\eps}^T \!\!\! \sup_{s\in [t-\eps, t]} |X_s - X_t|^2 dt \Big]}\\
		&\le&
		C \sqrt{\rho(\eps) + \eps},
	\end{eqnarray*}
	together with
	$\alpha_s \equiv 0$ for $s \ge T$,
	it follows from \eqref{eq:Regularity_gamma} that, for $\eps >0$,
	\begin{eqnarray*}
		\limsup_{n \to \infty} \E^{\P^n} \Big[\! \int_0^T \!\!\! \gamma_1(t, X) \alpha_t dt \Big]
		\!\!\!&\le&\!\!\!
		\limsup_{n \to \infty} 
		\Big(
		\E^{\P^n} \Big[
			\int_{\eps}^T \frac{1}{\eps} \int_{t-\eps}^t \gamma_1(s, X) ds \alpha_t dt
		\Big]
		+ C \delta(n ,\eps)
		\Big) \\
		\!\!\!&\le&\!\!\!
		\limsup_{n \to \infty} 
		\Big(
		\E^{\P^n} \Big[ \int_0^T \gamma_1(t, X) \frac{1}{\eps} \int_t^{t+ \eps} \alpha_s ds dt \Big] + C \delta(n ,\eps) \Big]
		\Big)\\
		\!\!\!&\le&\!\!\!
		\limsup_{n \to \infty} 
		\Big(
		\E^{\P^n} \Big[ \int_0^T \gamma_1(t, X)  \sqrt{ \frac{1}{\eps} \int_t^{t+ \eps} \alpha^2_s ds} dt \Big] + C \delta(n ,\eps)
		\Big)\\
		\!\!\!&=&\!\!\!
		\E^{\Ph} \Big[ \int_0^T \gamma_1(t, X)  \sqrt{ \frac{1}{\eps} \int_t^{t+ \eps} \alpha^2_s ds} dt \Big] + C \sqrt{\rho(\eps) + \eps} ,
	\end{eqnarray*}
	where the last equality is due to \eqref{eq:cvg_Pn} and the fact that $\sup_{n \ge 1} \E^{\P^n} [ \int_0^T \alpha^2_s ds] < \infty$.
	Letting $\eps \to 0$, we deduce that
	$$
		\limsup_{n \to \infty} \E^{\P^n} \Big[\! \int_0^T \!\!\! \gamma_1(t, X) \alpha_t dt \Big]
		~\le~
		 \E^{\Ph} \Big[\! \int_0^T \!\!\! \gamma_1(t, X) \alpha_t dt \Big].
	$$
	Together with \eqref{eq:cvg1} and \eqref{eq:cvg2}, this leads to 
	$$
		\E^{\Ph} \Big[ \Phi(X) - \int_0^T G_t(X, \alpha_t) dt \Big]
		~\ge~
		\limsup_{n \to \infty}
		\E^{\P^n} \Big[ \Phi(X) - \int_0^T G_t(X, \alpha_t) dt \Big]
		~=~
		\vr(0, x_0),
	$$
	and therefore $\Ph$ is an optimal control rule for problem \eqref{eq:dual_pb_equiv}.
	
	\vspace{0.5em}
	
	It remains to prove the claim in \eqref{eq:claim_tight}.
	Assume that \eqref{eq:claim_tight} is not true. Then there exist $M>0$, a sequence of positive constants $\theta_{n}\to 0$, 
	together with a sequence of stopping times $(\sigma_n, \tau_{n})_{n} \subset \Tc_M \x \Tc_M$ such that $\sigma_{n}\le \tau_{n}\le \sigma_{n}+\theta_{n}$ and 
	\begin{align}\label{eq: contra tight}
		2 c ~:=~ \liminf_{n}\E^{\P^{n}}[\int_{\sigma_{n}}^{\tau_{n}}|\alpha_{s}|^{2}ds] ~>~ 0.
	\end{align}
	Set
	$$
		\phi:=\vr-\Gammab_{\eps_0}
		~~\mbox{and}~~
		\xi_{n}:= \E^{\P^{n}}_{\sigma_{n}} \big[ \phi(\tau_{n},X)
		-
		\phi(\tau_{n}, (X \oplus_{\sigma_{n}}(X_{\tau_{n}}-X_{\sigma_{n}}))_{ \sigma_{n} \wedge \cdot }) \big].
	$$
	Using the same arguments as in the second part of the proof of Proposition \ref{prop : bar v - Gamma concave}, one has  $\E^{\P^{n}}[\xi_{n}]\to 0$ as $\theta_n \to 0$.
	Further, by \eqref{eq:non_increas} and the fact that $\phi$ is Dupire-concave, one obtains
	\begin{eqnarray*}
		&&
		\E^{\P^{n}}_{\sigma_{n}} \Big[ \vr(\tau_{n},X)-\frac12\int_{\sigma_{n}}^{\tau_{n}}\gamma_2(s, X_{s}) \alpha^{2}_{s}ds \Big] \\
		&=&\E^{\P^{n}}_{\sigma_{n}} \Big[\phi(\tau_{n},(X \oplus_{\sigma_{n}}(X_{\tau_{n}}-X_{\sigma_{n}}))_{\sigma_{n} \wedge \cdot})-\frac12\int_{\sigma_{n}}^{\tau_{n}}\eps_0 \alpha^{2}_{s}ds \Big]
		+\Gammab_{\eps_0}(\sigma_{n},X) + \xi_{n}\\
		&\le& \phi(\sigma_{n},X) + C\theta_{n} - \frac{\eps_0}2 \E^{\P^{n}}_{\sigma_{n}} \Big[\int_{\sigma_{n}}^{\tau_{n}} \alpha^{2}_{s}ds \Big] + \Gammab_{\eps_0}(\sigma_{n},X)+ \xi_{n}\\
		&=& \vr(\sigma_{n},X)+C\theta_n-\frac{\eps_0}2\E^{\P^{n}}_{\sigma_{n}} \Big[ \int_{\sigma_{n}}^{\tau_{n}} \alpha^{2}_{s}ds \Big]+\xi_{n}.
	\end{eqnarray*}
	Now recall that $\E^{\P_{n}}[\int_{\sigma_{n}}^{\tau_{n}} |\alpha_{s}|ds]\le C\theta_{n}^{\frac12}$ by \eqref{eq:bounded_Phi_X2}, for some $C>0$. Thus, \eqref{eq: contra tight} leads to
	$$
		\E^{\P^{n}} \Big[ \vr(\tau_{n},X) - \int_{\sigma_{n}}^{\tau_{n}} G_s(X, \alpha_s) ds \Big]
		~\le~ 
		\E^{\P^{n}}[\vr(\sigma_{n},X)]- \eps_0 c  +O(\theta_{n}^{1/2})+\xi_{n},
	$$
	for $n$ large enough, and therefore 
	\begin{equation} \label{eq:lim_vr}
		\lim_{n \to \infty}
		\E^{\P^{n}} \Big[ \vr(\tau_{n},X) - \int_{\sigma_{n}}^{\tau_{n}} G_s(X, \alpha_s) ds \Big]
		~\le~
		\lim_{n \to \infty}
		\E^{\P^n} \big[ \vr(\sigma_{n},X) \big] - \eps_0 c.
	\end{equation}
	On the other hand, it follows from the dynamic programming principle in Proposition \ref{prop:equiv_dpp} that 
	\begin{equation}
		\vr(0, x_0)
		\ge
		\lim_{n}\E^{\P^{n}} \Big[ \vr(\tau_n,X) \!-\!\! \int_{0}^{\tau_n}\!\! G_s(X, \alpha_s)ds \Big]
		\ge
		\lim_{n}\E^{\P^{n}} \Big[ \Phi(X) \!-\!\! \int_{0}^T\!\! G_s(X, \alpha_s)ds \Big]
		= \vr(0,x_0).
	\end{equation}	
	The above implies that 
	$$
		\lim_{n \to \infty}
		\E^{\P^{n}} \Big[ \vr(\tau_{n},X) - \int_{\sigma_{n}}^{\tau_{n}} G_s(X, \alpha_s) ds \Big]
		~=~
		\lim_{n \to \infty}
		\E^{\P^n} \big[ \vr(\sigma_{n},X) \big],
	$$
	which is a contradiction to \eqref{eq:lim_vr}.
	Therefore, the claim in \eqref{eq:claim_tight} holds true.
	\endproof

	\begin{remark}[Sufficient conditions for \eqref{eq: hyp eps conc}]\label{rem: eps conca}  
	Consider the Markovian setting in which $\Phi(\xr)=\Phi(\xr_{T}\1_{\{T\}})=: \bar \Phi(\xr_{T})$ and $G(t,\xr,a)=G(t,\xr_{t}\1_{\{t\}},a)=:\bar G(t,\xr_{t},a)$ for $t\le T$ and $\xr \in D([0,T])$, $a\in \R$. Then, the Hamilton-Jacobi-Bellman equation associated to $\bar \vr$, defined by $\bar \vr(t,\cdot):=\vr(t,\cdot\1_{\{t\}})$ for $t\le T$, is 
	\begin{align}\label{eq: HJB}
	\partial_{t}\vp+\bar F(\cdot,\partial^{2}_{xx}\vp)=0,\;\;\mbox{on } [0,T)\x \R,
	\end{align}
	with terminal condition $\vp(T,\cdot)=\bar \Phi$, in which 
	\begin{align}\label{eq: Fenchel}
	\bar F(t,x,z):=\sup_{a\in \R}\frac12\left(a^{2}z-G(t,x,a)\right),\;(t,x,z)\in [0,T]\x \R^{2}.
	\end{align}
	Assume further that $(t,x)\in [0,T]\x \R \mapsto \gamma_{2}(t,x\1_{\{t\}})$ is uniformly continuous, bounded, with bounded inverse, that  $\bar F$ is uniformly continuous on each $\Dc_{n}:=\{\bar F\le n\}$, $n\ge 1$, with $(\partial_{t} \bar F/\bar F)^{-}$  bounded on $\Dc:=\cup_{n\ge 1}\Dc_{n}$, that $\bar F(\cdot,0)=0$, $\bar F\in \cap_{n\ge 1} C^{1,3,3}_{b}(\Dc_{n}^{n})$ with $\Dc_{n}^{n}:=\Dc_{n}\cap ([0,T]\x \R\x [-n,n]])$, $\partial_{z}\bar F>0$ on $\Dc$, that  $|\partial_{z} \bar F|+|\partial_{z} \bar F|^{-1}$ is bounded on each $\Dc_{n}^{n}$, $n\ge 1$, and that $z\in (-\infty,0]\mapsto \partial_{x} \bar F(\cdot,z)$ has at most linear growth (uniformly in the other variables). Then, the a-priori estimates of \cite[Proposition 3.10]{BLSZ18} are valid, as well as the arguments in the proof of \cite[Theorem 3.11]{BLSZ18} when
	
	\;~~~\;{\rm (i)} $\bar \Phi\in C^{2+\iota}_{b}$ for some $\iota>0$, 
	 
	\;~~~\;{\rm (ii)} $\{(t,x,\partial^{2}_{xx}\bar \Phi(x)): (t,x)\in [0,T]\x\R\}\subset \Dc_{n}$, for some $n\ge 1$.
	 \vs2
	 
	Namely, in this case, \eqref{eq: HJB} admits a unique solution $u \in C^{1,2}_{b}([0,T]\x \R)$ such that the sup in the definition of $\bar F(t,x,\partial^{2}_{xx} u(t,x))$ is achieved and bounded uniformly in $(t,x)\in [0,T]\x \R$. By immediate verification, we have $u=\bar \vr$ and a strong solution to our dual optimal control problem exists. 
 		Moreover,  by \cite[Theorem 3.11]{BLSZ18} again, $u$ satisfies \eqref{eq: hyp eps conc} for some $\eps_{0}>0$ that only depends on $ \eps_{\bar \Phi}>0$ such that $x\mapsto \bar \Phi(x)-\bar \Gamma_{\eps_{\bar \Phi}}(T,x)$ is concave.
 
	 \vs2
	 
	Assume now that one can find a sequence $(\bar \Phi_{n},\bar G_{n})_{n\ge 1}$   whose elements each  satisfy the above requirement and converge  pointwise in a monotone way to $(\bar \Phi,\bar G)$ with $\bar \Phi_{n}\le \bar \Phi $, $\bar G_{n}\ge \bar G $ and such that one can choose $\eps_{\bar \Phi_{n}} = \eps_{\bar \Phi}>0$. Then, it is not difficult to see that the corresponding sequence of solutions $(u_{n})_{n\ge 1}$ converges pointwise to $\bar \vr$. From the preceding discussion, it follows that $\bar \vr$ satisfies \eqref{eq: hyp eps conc} as well with $\eps_{0}=\eps_{\bar \Phi}$ (although $\bar \vr$ will not be smooth in general and a strong solution to \eqref{eq: HJB} may not exist). We refer to the proofs of \cite[Theorem 3.11 and Lemma 3.13]{BLSZ18} for the construction of such approximating sequences.
	
	\vs2

	In general, we might have to face-lift $\bar \Phi$ so as to ensure that the corresponding $\eps_{\bar \Phi}$ is positive. This can be done as explained in Remark \ref{rem: conc enveloppe pour cond term} with $\bar \Gamma_{\eps}$ in place of $\bar \Gamma_{0}$ for a small $\eps>0$. In this case, the hedging strategy will hedge the face-lifted payoff, and will (in general) only be a super-hedging strategy for the original one.

\vs2

	Extensions of the above arguments to path dependent situations in which the coefficients depend on $\xr$ also through quantities of the form $(\int_{0}^{t} \xr_{s}\rho_{s}ds)_{t\le T}$, for some continuous deterministic process $\rho$,  are straightforward.

	\end{remark}
 
	\begin{remark}[Comparison with \cite{BLSZ18}] 
		Conditions in \cite{BLSZ18} under which a solution to the primal hedging problem exists are far more restrictive than the ones of the previous Remark \ref{rem: eps conca}. This is due to the fact that \cite{BLSZ18} looks directly to a solution to the primal problem and requires a $C^{1,4}_{b}$-solution to \eqref{eq: HJB} to construct it (see \cite[Corollary 3.12]{BLSZ18}). Since we only want to establish  \eqref{eq: hyp eps conc}, we only need  here  solutions to approximating   dual problems, each associated to a $C^{1,2}_{b}$-value function,   but we do not care about the regularity of the limit.
	\end{remark}

 
\appendix

\section{Appendix: A version of It\^{o}'s formula for path-dependent functionals}

	We provide here a particular version of the It\^{o}'s formula for path-dependent functional as initially introduced by Dupire \cite{dupireito}, and then studied by \cite{cont2013functional}.
	 {See also}  \cite{saporito2018functional} for convex functionals. 
	Our versions exploit the fact that the functional $\vp: \Theta \to \R$ we consider is non-increasing in time, Dupire-concave, 
	and  produce a martingale along the path of a given martingale, up to some correction term. 
	This allows us to assume less regularity on the function, that is, $\vp $ {is only} $\Cb^{0,1}(\Theta)$  {in the spirit of \cite{coquet2006natural,russo1996ito} who considered  Markovian settings.}  
	For simplicity, the results will only be stated for a one dimensional process, but it clearly holds for $d$-dimensional ones.
 {We first assume that the functional itself is Dupire-concave and non-increasing in time, before to generalize this result in Proposition \ref{prop: Dupire Tanaka corrige} to match with the conditions under which we need to apply it in Proposition \ref{prop : ito surbar v} above.}
	\begin{proposition}\label{prop: Dupire Tanaka} 
		Let $(\Omega,\Fc,\F,\P)$ be a filtered probability, equipped with a continuous semi-martingale $Z$ such that 
		$\E^{\P} \big[ \|Z\|^2 \big] <\infty$.  
		Let $\vp: \TD\mapsto \R$ be a non-anticipative function in $ \Cb^{0,1}(\Theta)$ which is Dupire-concave and 
		such that $\nabla_{\xr} \vp$ is locally bounded. Assume further that  the map $s \mapsto \vp_s(\xr_{t \wedge \cdot})$ is non-increasing on $[t,T]$ for all $(t,\xr) \in \Theta$.
		Then, there exists a predictable non-increasing process $K$ starting at $0$ such that 
		$$
			\vp_{\cdot}(Z)=\vp_{0}(Z_{0})+\int_{0}^{\cdot} \nabla_{\xr}\vp_{t}(Z)dZ_{t}+K\;,\;\mbox{ on }[0,T]. 
		$$
	\end{proposition}
	\proof Let us adapt the proof of \cite[Theorem 4.1]{cont2013functional} to our context. 
	Fix $t^{n}_{i} = ih^{n}$ for $i= 0, 1, \cdots, n$ with $h^{n}:=T/n$, $n\ge 1$.
	Set $Z^{n} := \sum_{i} Z_{ t^{n}_{i}}\1_{[t^{n}_{i},t^{n}_{i+1})}$.   
	Then,   
	\begin{align*} 
		\vp_{t^{n}_{i+1}}(Z^{n})-\vp_{t^{n}_{i}}(Z^{n})=  \vp_{t^{n}_{i+1}}(Z^{n})-\vp_{t^{n}_{i+1}}(Z^{n}_{\wedge t^{n}_{i}})+
		\vp_{t^{n}_{i+1}}(Z^{n}_{\wedge t^{n}_{i}})-\vp_{t^{n}_{i}}(Z^{n}).
	\end{align*}
	Applying the Meyer-Tanaka formula to  $r\in [t^{n}_{i},t^{n}_{i+1}]\mapsto \vp_{t^{n}_{i+1}}(Z^{n}_{\wedge t^{n}_{i}}\oplus_{t^{n}_{i+1}}(Z_{r}-Z_{t^{n}_{i}}))$, 
	it follows that
	\begin{align*}
		&\vp_{t^{n}_{i+1}}(Z^{n})-\vp_{t^{n}_{i+1}}(Z^{n}_{\wedge t^{n}_{i}})\\
		&=\int_{t^{n}_{i}}^{t^{n}_{i+1}} \nabla_{\xr} \vp_{t^{n}_{i+1}}(Z^{n}_{\wedge t^{n}_{i}}\oplus_{t^{n}_{i+1}}(Z_{r}-Z_{t^{n}_{i}})) dZ_{r} +K^{n}_{t^{n}_{i+1}}-K^{n}_{t^{n}_{i}}
	\end{align*}
	in which $K^{n}$ is a non-increasing predictable process starting at $0$. 
	On the other hand,   $\vp_{t^{n}_{i+1}}(Z^{n}_{\wedge t^{n}_{i}}) -\vp_{t^{n}_{i}}(Z^{n})\le 0$ by our monotonicity assumption.
	Hence, for $t\le T$,
	\begin{align*}
		\vp_{t}(Z)=& \vp_{0}(Z) + \sum_{j=0}^{n}\1_{[t^{n}_{j},t^{n}_{j+1})}(t)\left(\vp_{t}(Z) - \vp_{t^{n}_{i}}(Z^{n}) \right)  \\
		&+ \sum_{j=0}^{n}\1_{[t^{n}_{j},t^{n}_{j+1})}(t)\sum_{i=0}^{(j-1) \wedge (n-1)} \int_{t^{n}_{i}}^{t^{n}_{i+1}} \nabla_{\xr} \vp_{t^{n}_{i+1}}(Z^{n}_{\wedge t^{n}_{i}}\oplus_{t^{n}_{i+1}}(Z_{r}-Z_{t^{n}_{i}})) dZ_{r} 
		+ \widetilde K^{n}_{t},
	\end{align*}
	in which $\widetilde K^{n}$ is another non-increasing predictable process starting at $0$. 
	By \cite[Lemma A3]{cont2013functional}, $\Ninfty{Z^{n}-Z}\to 0$ as $n\to \infty$, so that  
	$$
		\sum_{j=0}^{n}\1_{[t^{n}_{j},t^{n}_{j+1})}(t)\left(\vp_{t}(Z) - \vp_{t^{n}_{i}}(Z^{n})\right)\to 0\;\;\mbox{ as } n\to \infty,\;\P\mbox{-a.s.}
	$$
 	Moreover, up to an additional localization argument, one can assume that 
 	$ \nabla_{\xr} \vp$ is bounded. Thus, the above, \cite[Chapter IV, Theorem 32]{protter} and the fact that $\vp\in \Cb^{0,1}(\Theta)$ imply that 
 	$$
 		\Big\| 
			\sum_{j=0}^{n}\1_{[t^{n}_{j},t^{n}_{j+1})}(\cdot)  \sum_{i=0}^{(j-1) \wedge (n-1)} \int_{t^{n}_{i}}^{t^{n}_{i+1}} \nabla_{\xr} \vp_{t^{n}_{i+1}}(Z^{n}_{\wedge t^{n}_{i}}\oplus_{t^{n}_{i+1}}(Z_{r}-Z_{t^{n}_{i}})) dZ_{r}  
			- \int_{0}^{\cdot}  \nabla_{\xr} \vp_{r}(Z) dZ_{r} 
		\Big\|
		\to 0  
	$$
	in probability  as $n\to \infty$, and therefore $\P$-a.s. along a subsequence. 
	Finally, $(\tilde K^{n})_{n\ge 1}$ being a sequence of non-increasing processes starting at $0$, it converges to some  non-increasing predictable process $K$ starting at $0$. 
	Combining the above, we obtain that
	\begin{align*}
		\vp_{\cdot}(Z) 
		=
		\vp_{0}(Z)  +  \int_{0}^{\cdot} \nabla_{\xr} \vp_{r}(Z) dZ_{r}   + K.
	\end{align*}
	\endproof
	
	\vspace{0.5em}

	We now state a variant result of Proposition \ref{prop: Dupire Tanaka},
	in which we alleviate the conditions of monotonicity and concavity. 
	Since the proof is very similar to the one of Proposition \ref{prop: Dupire Tanaka}, we only explain the main modifications.

	\begin{proposition}\label{prop: Dupire Tanaka corrige} 
		Let $(\Omega,\Fc,\F,\P)$ be a filtered probability, equipped with a continuous semi-martingale $Z$ such that  $\E^{\P}[\Ninfty{Z}^{2}]<\infty$.  
		Let $\vp: \Theta \mapsto \R$ be a non-anticipative map in $ \Cb^{0,1}(\Theta)$ such that $\nabla_{\xr} \vp$ is locally bounded.
		Assume that there exists $R \in \Cb^{1,2}(\Theta)$ and a continuous function $\ell: [0,T]\to \R$ such that:
		\begin{enumerate}[\rm (1)]
			\item $\vp-R$ is Dupire-concave.
			\item $s \mapsto \vp_s(\xr_{t \wedge \cdot}) -\ell(s)$ is non-increasing on $[t,T]$,  for any $(t, \xr) \in \Theta$. 
		\end{enumerate} 
		Then, there exists a non-increasing predictable process $K$ starting at $0$ such that 
		$$
			\vp_{\cdot}(Z)-\int_{0}^{\cdot }\frac12 \nabla^{2}_{\xr} R_{r}(Z)d\langle Z\rangle_{r}=\vp_{0}(Z)+  \int_{0}^{\cdot} \nabla_{\xr} \vp_{r} (Z)dZ_{r} + K +\ell(\cdot)-\ell(0).
		$$
		Moreover, if  $Z$ and  $\vp_{\cdot}( Z) -B$ are $(\P,\F)$-martingales, for some predictable bounded variation process $B$, then 
		$$
			\vp_{\cdot}(Z)=\vp_{0}(Z_{0})+\int_{0}^{\cdot} \nabla_{\xr}\vp_{t}(Z)dZ_{t}+B\;,\;\mbox{ on }[0,T]. 
		$$
	\end{proposition}
	\proof 
	By a straightforward adaptation of the arguments of the proof of Proposition \ref{prop: Dupire Tanaka}, one can find a non-increasing predictable process $K$, starting at $0$, such that 
	\begin{align*}
		(\vp_{\cdot}-R_{\cdot})(Z)-\ell(\cdot) 
		= (\vp_{0}-R_{0})(Z)-\ell(0)  +  \int_{0}^{\cdot} (\nabla_{\xr} \vp_{t}-\nabla_{\xr}R_{t})(Z) dZ_{t}   + K-\int_{0}^{\cdot} \partial_{t}R_{t}(Z)dr.
	\end{align*}
	By the It\^{o}'s formula for smooth path-dependent functionals, \cite[Theorem 4.1]{cont2013functional}, 
	$$
		R_{\cdot}(Z)=R_{0}(Z)+\int_{0}^{\cdot }\nabla_{\xr}R_{t}(Z) dZ_{t}+\int_{0}^{\cdot }\frac12 \nabla^{2}_{\xr} R_{t}(Z)d\langle Z\rangle_{t} +\int_{0}^{\cdot }\partial_{t}R_{t}(Z) dr.
	$$
	Hence, 
	$$
		\vp_{\cdot}(Z)=\vp_{0}(Z)+  \int_{0}^{\cdot} \nabla_{\xr} \vp_{t} (Z)dZ_{t} + \int_{0}^{\cdot }\frac12 \nabla^{2}_{\xr} R_{t}(Z)d\langle Z\rangle_{t}+ K +\ell(\cdot)-\ell(0).
	$$
	If $Z$ and $\vp_{\cdot}(Z)-B$ are martingales, for some predictable bounded variation process $B$, then $ \int_{0}^{\cdot }\frac12 \nabla^{2}_{\xr} R_{t}(Z)d\langle Z\rangle_{t}+K +\ell(\cdot)-\ell(0)\equiv B$.
	\endproof

 \def\cprime{$'$} \def\cprime{$'$}

 \end{document}